# A Moving Discontinuous Galerkin Finite Element Method with Interface Conservation Enforcement for Compressible Flows

Hong Luo[1] and Gianni Absillis[2]
*North Carolina State University, Raleigh, NC 27695, USA*

Robert Nourgaliev[3]
*Lawrence Livermore National Laboratory, Livermore, CA 94550, USA*

## Abstract

A moving discontinuous Galerkin finite element method with interface conservation enforcement (MDG+ICE) is developed for solving the compressible Euler equations. The MDG+ICE method is based on the space-time DG formulation, where both flow field and grid geometry are considered as independent variables and the conservation laws are enforced both on discrete elements and element interfaces. The element conservation laws are solved in the standard discontinuous solution space to determine conservative quantities, while the interface conservation is enforced using a variational formulation in a continuous space to determine discrete grid geometry. The resulting over-determined system of nonlinear equations arising from the MDG+ICE formulation can then be solved in a least-squares sense, leading to an unconstrained nonlinear least-squares problem that is regularized and solved by Levenberg-Marquardt method. A number of numerical experiments for compressible flows are conducted to assess the accuracy and robustness of the MDG+ICE method. Numerical results obtained indicate that the MDG+ICE method is able to implicitly detect and track all types of discontinuities via interface conservation enforcement and satisfy the conservation law on both elements and interfaces via grid movement and grid management, demonstrating that an exponential rate of convergence for Sod and Lax-Harden shock tube problems can be achieved and highly accurate solutions without overheating to both double-rarefaction wave and Noh problems can be obtained.

## I. Introduction

The discontinuous Galerkin (DG) finite element methods[1-30] have recently become a popular choice to solve conservation laws with arbitrary order of accuracy in Eulerian, Lagrangian, and Arbitrary Lagrangian-Eulerian (ALE) formulations. Originally introduced for the solution of neutron transport equations[1], nowadays they are widely used in computational fluid dynamics (CFD), computational acoustics, and computational magneto-hydrodynamics. The discontinuous Galerkin methods have many attractive features: (1) They have several useful mathematical properties with respect to conservation, stability, and convergence; (2) The methods can be easily extended to higher-order (>2nd) approximations; (3) The methods are well suited for complex geometries since they can be applied on unstructured grids. In addition, the methods can also handle non-conforming elements, where the grids are allowed to have hanging nodes; (4) The methods are highly parallelizable, as they are compact and each element is independent; (5) They can easily handle adaptive strategies, since refining or coarsening a grid can be achieved without considering the continuity restriction commonly associated with the conforming elements. The methods allow for easy implementation of *hp*-refinement, for example, where the order of accuracy, or shape, can vary from element to element; (6) They have the ability to compute low Mach number flow problems without recourse to the time-preconditioning techniques normally required by the finite volume methods; (7) they provide the powerful feature of adjoint consistency for adjoint-based optimization. Furthermore, the space-time discontinuous Galerkin methods[31-34], which do not separate space and time, simply regard a time-dependent n-dimensional problem as a n+1-dimensional steady problem and use discontinuous basis functions across element faces, both in space and time in the finite element approximation. In addition to inheriting all the features from the standard DG methods, the space-time DG finite element methods can provide uniform or different arbitrary order of accuracy in both space and time and are naturally suited for moving/deforming boundary problems, as the geometric conservation law is automatically satisfied and the issue of maintaining conservative interpolation does not exist, when a remeshing occurs or the topology of a mesh changes.

---

[1] Professor, Department of Mechanical and Aerospace Engineering
[2] Ph.D. student, Department of Mechanical and Aerospace Engineering
[3] Computational Physicist

However, one area where the DG methods struggle to deliver the designed optimal order of accuracy and have to face the same challenging issue all other high-order methods face is for flows that are not smooth and contain discontinuous interfaces, such as shock waves, contact discontinuities, and material interfaces. Such flow features are ubiquitous in many science and engineering applications, particularly in high-speed flow and multi-material flow problems. A flow past a vehicle at hypersonic speeds is one example, where the computed surface heat flux exhibits extreme sensitivity to the resolution of a strong shock wave ahead of the vehicle. A common practice to address discontinuities is to use the so-called shock capturing approach. Schemes of this type are widely used on modern compressible CFD solvers, because they are algorithmically simple and easy to implement, perform well on current computer architectures, and provide reasonable accuracy for a wide range of problems. However, they have a fatal flaw in that the interface conservation is not considered and can never be achieved for discontinuities. Neglecting the interface conservation is the root cause of the two long outstanding issues for solving conservation laws: how to compute numerical fluxes at interfaces and how to avoid and suppress spurious oscillations in the vicinity of discontinuities for high-order methods (>1st). Continuous efforts in the last few decades for addressing these two issues manifest that the 'perfect' or 'right' answers to these two questions remain evasive, indicating a pressing need to explore new and innovative approaches that address these two issues or, even better, eliminate them entirely.

Moving grid methods[35-44] are widely used for solving a variety of flow problems in computational fluid dynamics for different reasons. Arbitrary Lagrangian-Eulerian (ALE) methods are developed for solving moving and deforming boundary problems[36]. In $r$-adaptation methods[37-38], mesh points are relocated into regions in need of high resolution in order to significantly increase solution accuracy. In the Lagrangian formulation, meshes are moved with the fluid velocity to explicit track contact discontinuities and material interfaces[29,35]. In shock-fitting or shock-tracking methods, grids are moved such that inter-element boundaries are alighted with shock waves via Rankine-Hugoniot or other conditions.[39-44]. In general, it is difficult to track these discontinuities, when their topologies are unknown a priori or when their topologies change in time. In addition, these interface-tracking methods are only designed to track or fit one particular type of discontinuity and therefore cannot track different types of interfaces simultaneously. For example, the Lagrangian formulation only tracks material interfaces, while shock fitting methods only detect and track shock waves. Moreover, all grid movements in all these applications should be considered as explicit and passive in the sense that they are not independent but rather explicitly driven by different and specific considerations, therefore making them algorithmically complex and difficult to implement.

Recently, a moving discontinuous Galerkin finite element method with interface condition enforcement, termed MDG-ICE, was developed by Corrigan et al. for compressible flows with discontinuous interfaces[45-48]. In the MDG-ICE formulation, both the flow field and discrete geometry are considered as independent variables. A space-time DG formulation is used to solve the governing equations in the standard discontinuous solution space, and the discrete geometry variables are determined by enforcing the interface condition in its discontinuous solution trace space. Consequently, all types of discontinuities are implicitly detected and tracked by making interface fitting an intrinsic part of the underlying MDG-ICE formulation. Since the interface condition is enforced in the MDG-ICE formulation, the space-time flux function across an element interface is uniquely defined. Consequently, the MDG-ICE method does not need to use a common numerical flux function in the form of an exact or approximate Riemann solution in order to maintain its linear stability. Furthermore, since all types of discontinuities are implicitly tracked and fitted, no strategies in the form of a limiter or ENO/WENO reconstruction, or an artificial viscosity, are required to eliminate spurious oscillations in the vicinity of discontinuities and thus maintain the nonlinear stability of the MDG method. The MDG-ICE method is especially advantageous and attractive in the sense that these two problems simply do not exist. Furthermore, a combination of the space-time DG formulation for solving the conservation laws and the interface conservation enforcement for governing mesh movement renders the MDG-ICE method a unique capability to implicitly track all types of discontinuities: shocks, contact discontinuities, and even locations where solution derivatives are discontinuous such as tails and heads of rarefaction waves. However, this MDG-ICE formulation leads to an over-determined system of nonlinear equations, which can be difficult and expensive to solve. A different MDG-ICE formulation leading to a determined system of nonlinear equations, designed to reduce the number of the nonlinear equations for both governing equations and interface condition, was introduced in Reference 49 to solve the scalar conservation laws. This formulation can be extended to the system of conservation laws, which will be presented in a follow-up paper.

The objective of the efforts presented in this work is to develop a moving discontinuous Galerkin finite element method with interface conservation enforcement (MDG+ICE) for solving the compressible Euler equations. Using the interface conservation enhancement instead of the interface condition enhancement is merely intended to better reflect and understand the idea, insight, and design of the MDG+ICE formulation, as the interface condition is simply the direct result of enforcing the conservation laws on interfaces. In this MDG+ICE formulation, the compressible Euler equations are solved using a standard space-time DG formulation, while the geometric variables are determined by enforcing the interface conservation using a continuous variational formulation. As a result, the number of nonlinear equations for the interface condition is significantly reduced. The resulting over-determined system of nonlinear equations arising from the MDG+ICE formulation is then solved in a least-squares sense, leading to an unconstrained nonlinear least-squares problem, which is regularized, solved by Levenberg-Marquardt method. A number of numerical experiments for both 1D unsteady and 2D steady compressible flow problems are conducted to assess accuracy and robustness of the developed MDG method. Numerical results obtained indicate that the MDG+ICE method is able to deliver the designed optimal order of convergence even for discontinuous solutions and solutions with singularities, detect all types of discontinuities via interface conservation enforcement, and satisfy the compressible Euler equations and the interface conservation via grid movement and management. The superior accuracy of the MDG method is attributed to the enforcement of the interface conservation, while the robustness of the MDG method is attributed to the renunciation of using an upwind flux function to achieve the linear stability and limiters/ENO/WENO/Artificial viscosity strategies to achieve the nonlinear stability. The remainder of this paper is organized as follows. The governing equations are described in Section 2. The developed MDG+ICE method is presented in Section 3. Numerical experiments are reported in Section 4. Concluding remarks and future work are given in Section 5.

## II. Governing Equation

The compressible Euler equations governing unsteady inviscid compressible flows can be expressed in conservative form as

$$\frac{\partial \mathbf{U}}{\partial t} + \frac{\partial \mathbf{f}_i(\mathbf{U})}{\partial x_i} = 0, \tag{2.1}$$

where the summation convention is used. The conservative variable vector $\mathbf{U}$, and inviscid flux vector $\mathbf{f}$ are defined by

$$\mathbf{U} = \begin{pmatrix} \rho \\ \rho u_i \\ \rho e \end{pmatrix} \qquad \mathbf{f}_i = \begin{pmatrix} \rho u_i \\ \rho u_j u_i + p\delta_{ij} \\ u_i(\rho e + p) \end{pmatrix} \tag{2.2}$$

Here $\rho$, $p$, and $e$ denote the density, pressure, and specific total energy of the fluid, respectively, and $u_i$ is the velocity of the flow in the coordinate direction $x_i$. The pressure is determined by the equation of state

$$p = (\gamma - 1)\rho \left(e - \frac{1}{2}u_i u_i\right) \tag{2.3}$$

which is valid for a perfect gas, and where $\gamma$ is the ratio of the specific heats. We are mainly interested in solving the 1D unsteady and 2D steady compressible Euler equations. For the 1D unsteady compressible Euler equations, Eq. (2.1) can be written as a 2D problem in the *x-t* domain, $\Omega(x,t)$

$$\frac{\partial \mathbf{F}_x}{\partial x} + \frac{\partial \mathbf{F}_t}{\partial t} = 0, \tag{2.4}$$

where the two independent variables are time *t* and position *x*, and the *x*- and *t*-components of the flux vector $\mathbf{F}$ are

$$\mathbf{F}_x = \mathbf{f}, \text{and } \mathbf{F}_t = \mathbf{U}. \tag{2.5}$$

For the 2D steady-state compressible Euler equations, Eq. (2.1) can be expressed as a 2D problem in the *x-y* domain, $\Omega(x,y)$

$$\frac{\partial \mathbf{F}_x}{\partial x} + \frac{\partial \mathbf{F}_y}{\partial y} = 0, \tag{2.6}$$

where the two independent variables are position $(x,y)$, and the $x$- and $y$-components of the flux vector $\mathbf{F}$ are

$$\mathbf{F}_x = \mathbf{f}_x, \text{and } \mathbf{F}_y = \mathbf{f}_y. \tag{2.7}$$

### III. Moving Discontinuous Galerkin Method with Interface Conservation Enhancement

**3.1 Space-time discontinuous Galerkin formulation**

In the space-time discontinuous Galerkin formulation[31-34], the space and time are indistinguishable. A time-dependent problem can simply be regarded as a steady problem on a d+1-dimensional space-time mesh in a finite time interval [0,T], where d = 1, 2, 3 is the spatial dimension and +1 is the time dimension. Without loss of generality and for the sake of clarity, we present the space-time DG formulation for 1D unsteady Euler equations (2.4). We first introduce some notations. We assume that the space-time domain $\Omega(x,t)$ is subdivided into a collection of non-overlapping elements $\Omega_e(x,t)$. We use $\Gamma_e(x,t)$ to denote the boundary of $\Omega_e$ and $\mathbf{n}$ the unit outward normal vector to $\Gamma_e$. We introduce the following broken Sobolev space $\mathbf{V}_h^p$

$$\mathbf{V}_h^p = \{v_h \in [L_2(\Omega)]^m : v_h|_{\Omega_e} \in [\mathbf{V}_p]^m \; \forall \Omega_e \in \Omega \}, \tag{3.1}$$

which consists of discontinuous vector-valued polynomial functions of degree $p$, and where $m$ is the dimension of the unknown vector $\mathbf{U}$ and $\mathbf{V}_p$ is the space of all polynomials of degree $\leq p$. To formulate the discontinuous Galerkin method, we introduce the following weak formulation, which is obtained by multiplying the above conservation laws (2.4) by a test function $w_h$, integrating over an element $\Omega_e$, and then performing an integration by parts,

$$\begin{cases} \text{Find } \mathbf{U}_h \in \mathbf{V}_h^p \text{ such that} \\ -\int_{\Omega_e} \mathbf{F}(\mathbf{U}_h) \cdot \nabla w_h \; d\Omega + \int_{\Gamma_e} \mathbf{F}(\mathbf{U}_h) \cdot \mathbf{n} \, w_h \, d\Gamma = 0, \quad \forall \, w_h \in \mathbf{V}_h^p \end{cases} \tag{3.2}$$

where $\mathbf{U}_h$ and $w_h$ are represented by piecewise-polynomial functions of degree $p$, which are discontinuous between cell interfaces. Assume that $B_i$ is the basis of polynomial functions of degree $p$, Eq. (3.2) is then equivalent to the following system of N nonlinear equations,

$$\begin{cases} \text{Find } \mathbf{U}_h \in \mathbf{V}_p \text{ such that} \\ -\int_{\Omega_e} \mathbf{F}(\mathbf{U}_h) \cdot \nabla B_i \; d\Omega + \int_{\Gamma_e} \mathbf{F}(\mathbf{U}_h) \cdot \mathbf{n} \, B_i \, d\Gamma = 0, \quad 1 \leq i \leq N \end{cases} \tag{3.3}$$

where N is the dimension of the polynomial function space. Since the numerical solution $\mathbf{U}_h$ is discontinuous between element interfaces, the interface space-time fluxes are not uniquely defined, and need to be computed carefully for the consideration of stability. This scheme is called the discontinuous Galerkin method of degree $p$, or in short notation DG(P) method. By simply increasing the degree $p$ of the polynomials, the DG methods of corresponding higher order are obtained. The domain and boundary integrals in Eq. (3.3) are calculated using Gauss quadrature formulas. The number of quadrature points used is chosen to integrate exactly polynomials of order of $2p$ and $2p+1$ for the volume and surface inner products in the reference element. A numerical polynomial solution $\mathbf{U}_h$ in each element is expressed using a standard finite element basis as following

$$\mathbf{U}_h = \sum_{i=1}^{N} \mathbf{U}_i B_i(x,t) \tag{3.4}$$

where $B_i$ are a set of polynomial basis functions and $\mathbf{U}_i$ are the unknown coefficients to be determined. The space-time DG formulation inherits all features of the standard DG formulation, and in addition allows uniform or different order of accuracy in space and time, as the time and space are treated in the same way. The space-time DG formulation is especially attractive for the moving/deforming boundary problems, as the geometric conservation law is naturally satisfied. Furthermore, the issue of conservative interpolation does not exist when a remeshing occurs or the topology of a mesh is changed.

**3.2 Moving discontinuous Galerkin-Interface Conservation Enhancement Formulation**

Traditionally, the governing equations for the conservation laws are only solved on elements $\Omega_e$ (computational cells) as described in the previous sub-section. Similarly, the conservation laws can be enforced on element interfaces $\Gamma_e$.

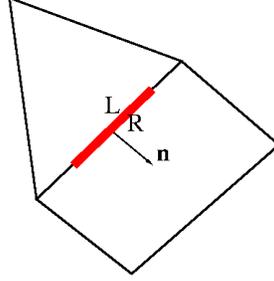

**Figure 1. Illustration of a zero-thickness control volume on an interface**

Applying the conservation laws on a zero-thickness control volume along an interface $\Gamma_e$ as shown in red in Figure 1 leads to the following jump condition for the flux function across the interface

$$\mathbf{F}(\mathbf{U}_h^R) \cdot \mathbf{n} - \mathbf{F}(\mathbf{U}_h^L) \cdot \mathbf{n} = \mathbf{0} \quad \Rightarrow \quad [\mathbf{F}(\mathbf{U}_h) \cdot \mathbf{n}] = \mathbf{0}, \tag{3.5}$$

where $\mathbf{U}_h^L$ and $\mathbf{U}_h^R$ are the conservative variable vector on the interface from the left and right elements respectively, and the bracket is the so-called the jump operator. This is the so-called interface condition, which is termed the transmission condition in the hybridized DG[50] or embedded DG[51] formulation. The interface conservation is never considered in all shock-capturing based schemes, because 1) it is automatically satisfied for smooth flows, as long as the flows are fully resolved, which can always be achieved by using an excessive mesh resolution; 2) it can never be satisfied for flows with discontinuities by simply using excessively refined meshes, as long as the discontinuities are not aligned with mesh interfaces. In other words, the interface conservation can only be achieved, if and only if the mesh interfaces are aligned with discontinuities. Recognizing the importance of the interface conservation enforcement and the necessity of a mesh movement to enforce the interface conservation, the main idea behind the MDG+ICE method is to treat the discrete geometry as an independent variable, which is determined by enforcing the interface conservation explicitly. Therefore, the MDG+ICE formulation is characterized by treating both flow field and discrete geometry as independent variables, and by solving the conservation laws on both elements and element interfaces simultaneously in the space-time domain. In the MDG+ICE method, discontinuous interfaces are not explicitly tracked and rather solved and obtained implicitly as a result of the interface conservation enforcement. Furthermore, what the MDG+ICE method offers is more than detecting and fitting discontinuities, which can also be achieved by other discontinuity fitting methods. By solving the conservation laws on interfaces, a continuous solution with a discontinuous derivative, such as the head and tail of a rarefaction wave, can be exactly resolved by the MDG method as demonstrated in next section, which cannot be obtained by other discontinuity tracking methods.

In the original MDG+ICE method introduced by Corrigan et al.[45-48], the geometry variables are determined by enforcing the interface condition Eq. (3.5) using a discontinuous variational formulation as follows,

$$\int_{\Gamma_e} [\mathbf{F}(\mathbf{U}_h) \cdot \mathbf{n}] \, v_h = \mathbf{0}, \tag{3.6}$$

where $v_h$ is a test function in the discontinuous solution trace space. The derivation of discrete equations for the interface condition is similar to the approach used in the so-called hybridized DG method[50]. The dimension of the discontinuous solution trace space is nface*nbasis, where nface is the number of faces in the grid, and nbasis is the number of basis functions for the DG(P) approximation on each face. Figure 2 illustrates the dimension of the test function space for the interface condition in the DG(P1) method, which is the number of red dots. In this case, nbasis is 2, as two basis functions are used to represent a linear polynomial in 1D. Consequently, the number of nonlinear equations derived from Eq. (3.6), i.e, the size of the nonlinear residual vector for the enforcement of the interface condition, is nface×nbasis×neqns, where neqns is the dimension of the flux vector **F**, i.e., d+2. (3, 4, and 5 for the 1, 2, and 3D Euler equations, respectively). The number of the geometric variables is 2×npoin, where npoin designates the number of grid points, i.e., the *x*- and *t*-coordinates of the grid vertices. Since the number of equations is larger than the number of unknowns, the veriational formulation (3.6) leads to an over-determined system of nonlinear equations. As both space and time are treated equally and the interface condition is enforced in the discontinuous solution trace space, the size of the resulting system of nonlinear equations can be large.

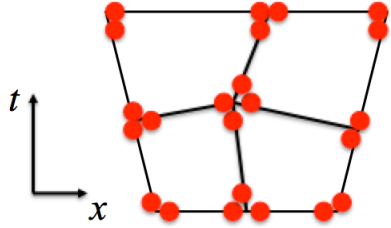

**Figure 2. Illustration for the dimension of the test function space for the interface condition in the case of a DG(Q1) approximation**

Alternative strategies can be explored and used to enforce the interface conservation. In this work, meshes are set uniform in time, i.e., cannot move in the *t*-direction. Furthermore, the geometric variables are determined by enforcing the interface conservation using a continuous variational formulation as follows,

$$\int_{\Gamma_e} [\mathbf{F}(u_h) \cdot \mathbf{n}]\, w_h = 0, \quad (3.7)$$

where $w_h$ is a test function in the continuous solution trace space. The derivation of discrete equations for the interface condition is similar to the approach used in the so-called embedded DG method[51]. Figure 2 illustrates the dimension of the test function space for the interface condition in the case of a CG(P1) approximation. In this case, the number of nonlinear equations obtained from Eq. (3.7) is npoin×neqns, where npoin is the number of grid points, i.e., the number of red dots in Figure 3. The number of the geometric variables is only npoin, i.e., the *x*-coordinate of the grid vertices, as the mesh cannot move in the *t*-direction. The number of equations derived for the interface condition is significantly reduced in this formulation. In the case of a scalar conservation law, the number of nonlinear equations for the interface condition is the same as the number of the geometric variables, therefore leading to a square system of nonlinear equations. Since meshes are uniform in time, the interface condition does not need to be enforced on the faces that separate the time slabs, i.e., all horizontal faces as shown in Figure 3. In this case, fluxes on these faces in Eq. (3.3) are computed in a naturally upwind manner, as the information can only be propagated from the past to the future according to the causality principle.

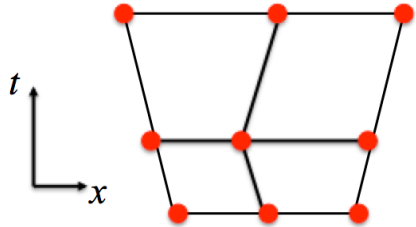

**Figure 3. Illustration for the dimension of the test function space for the interface condition in the case of a CG(Q1) approximation**

In our implementation, one element in time with a uniform time-step size is used, which leads to the traditional space-time DG formulation[31-33] normally described on one single space-time slab [$t_n$, $t_{n+1}$] as shown in Figure 4. Our justifications of using this MDG+ICE formulation are following: Firstly, the size of the resultant system of nonlinear equations is significantly reduced. For the DG(Q1) approximation, the number of unknown flow variables is ndegr×neqns×nelem, and the number of geometric variables is npoin/2, i.e., the number of the grid points at time step n+1, as illustrated in Figure 3, where ndegr designates the number of degrees of freedom in DG(P) approximation (3 for DG(P1) and 4 for DG(Q1) in general). An over-determined system of nonlinear equations for (ndegr×neqns×nelem+npoin/2) solution unknowns needs to be solved at each time step, thus significantly reducing both storage requirements and computational costs and ultimately resulting in a more efficient MDG+ICE method. Secondly, the implementation of the MDG method for the 3D unsteady problems becomes easy and straightforward. Otherwise, one would need to work in 4D space, which can be challenging and intimidating. Thirdly, this MDG+ICE formulation advances solutions using one cell in time and with a uniform time-step size, leading to a time-marching-like method, and consequently is easy to understand, analyze, implement, and solve. Indeed, we can prove that this MDG+ICE method has the property of resolving and tracking both contact and shock wave discontinuities exactly, which is given by the following two lemmas.

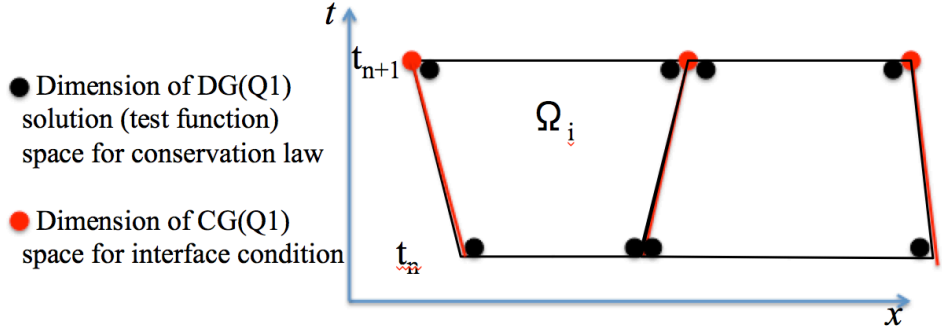

**Figure 4. Illustration for the dimension of the DG(Q1) solution for conservation law and the dimension of CG(Q1) for interface condition**

**Lemma 1**.  A contact discontinuity is exactly resolved (i.e., it moves with the convective velocity), if and only if the interface condition is satisfied

We prove lemma 1 by considering the following Riemann problem for the advection equation

$$\begin{cases} \dfrac{\partial u(x,t)}{\partial t} + a\dfrac{\partial u(x,t)}{\partial x} = 0 \\ u(x,t^n) = \begin{cases} u_l & \text{if } x < x_i \\ u_r & \text{if } x > x_i \end{cases} \end{cases} \tag{3.8}$$

where a single jump occurs at $x=x_i$.

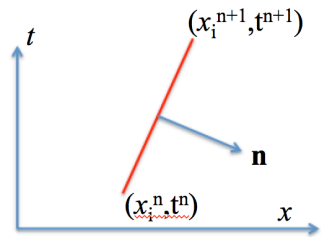

**Figure 5. Illustration for enforcing the interface condition along a discontinuity**

The lemma 1 is proved, if a contact discontinuity at $x_i$ and $t^n$ moves to a new position $x_i^{n+1} = x_i^n + a(t^{n+1} - t^n)$ at time $t^{n+1}$. The interface jump condition along a contact discontinuity as shown in Figure 4 is given by

$$[\mathbf{F}(u_h) \cdot \mathbf{n}] = [f]n_x + [u]n_t = a[u](t^{n+1} - t^n) + [u](x_i^{n+1} - x_i^n) = 0 \tag{3.9}$$

Therefore,

$$x_i^{n+1} = x_i^n + a(t^{n+1} - t^n), \tag{3.10}$$

which completes the proof of Lemma 1.

**Lemma 2**. A shock wave is exactly resolved (i.e., it moves with its velocity), if and only if the interface jump condition is satisfied.

We prove lemma 2 by considering the following Riemann problem for the inviscid Burger's equation

$$\begin{cases} \frac{\partial u(x,t)}{\partial t} + \frac{\partial}{\partial x}(\frac{u^2(x,t)}{2}) = 0 \\ u(x,t^n) = \begin{cases} u_l & \text{if } x < x_i \\ u_r & \text{if } x > x_i \end{cases} \end{cases} \quad (3.11)$$

where a single jump occurs at $x=x_i$.

The lemma 2 is proved, if a shock wave at $x_i$ and $t^n$ moves to a new position $x_i^{n+1} = x_i^n + \frac{u_r+u_l}{2}(t^{n+1} - t^n)$ at time $t^{n+1}$. The interface jump condition along the interface as shown in Figure 4 is given by

$$[\mathbf{F}(u_h) \cdot \mathbf{n}] = [f]n_x + [u]n_t = [\frac{u^2}{2}](t^{n+1} - t^n) + [u](x_i^{n+1} - x_i^n) = 0 \quad (3.12)$$

Therefore,

$$x_i^{n+1} = x_i^n + \frac{[\frac{u^2}{2}]}{[u]}(t^{n+1} - t^n) = x_i^n + \frac{u_r + u_l}{2}(t^{n+1} - t^n), \quad (3.13)$$

which completes the proof of Lemma 2.

### 3.3 Numerical space-time fluxes

Since the numerical solution $\mathbf{U}_h$ is discontinuous between element interfaces, the interface space-time fluxes in the surface integral in Eq. (3.3) are not uniquely defined, and therefore need to be computed carefully for the consideration of consistency and stability. In the standard discontinuous Galerkin and space-time discontinuous Galerkin formulation, the flux function is replaced by a common numerical flux based on an exact or approximate Riemann solution, which is absolutely required in order to be able to maintain the stability. However, in the MDG+ICE formulation, the interface space-time fluxes are uniquely defined by the enforcement of interface condition Eq. (3.5), which indicates

$$\mathbf{F}(\mathbf{U}_h^L) \cdot \mathbf{n} = \mathbf{F}(\mathbf{U}_h^R) \cdot \mathbf{n}. \quad (3.14)$$

Therefore, the flux function is not substituted by a common numerical flux function and instead the interior flux is retained. This means that the information exchange between two neighboring cells, i.e., the global coupling, is only achieved through the enforcement of the interface condition. Instead of using the separate flux function, a common flux function as the average of the left and right flux functions

$$\mathbf{F}(\mathbf{U}_h) \cdot \mathbf{n} = \frac{1}{2}(\mathbf{F}(\mathbf{U}_h^L) \cdot \mathbf{n} + \mathbf{F}(\mathbf{U}_h^R) \cdot \mathbf{n}) \quad (3.15)$$

is used in our MDG+ICE formulation in order to enhance the stability and coupling for under-resolved flows. This simple average flux function satisfies the so-called consistency requirement for the numerical fluxes in the MDG+ICE formulation, i.e., Eq. (3.14) and is found to enhance the stability of the MDG+ICE method. It should be emphasized that the MDG+ICE formulation does not need an upwind flux function from a Riemann solution, as solving the conservation laws on interfaces explicitly enforces the interface conservation. In fact, how to extend an upwind flux to the space-time DG formulation, that satisfies the consistence requirement (3.14), remains an open problem, although most of the upwind flux functions can be directly used in the steady-state MDG+ICE formulation.

### 3.4 Solving the nonlinear least-squares problem

As can be seen, the MDG+ICE formulation leads to an over-determined system of nonlinear equations that needs to be solved at each time step or more precisely for one space-time slab

$$\mathbf{R}(\mathbf{U}) = \mathbf{0}, \tag{3.16}$$

where **U** is the unknown solution vector, which includes both flow variables and geometric variables, and **R** represents the nonlinear residual function, which includes both residual from the space-time DG formulation and residual from the enforcement of interface condition. Let the dimension of the solution vector **U** be n and the dimension of the nonlinear residual vector **R** be m. Since m>n, the over-determined system of nonlinear equations (3.16) is solved in a least-squares sense

$$\min \tfrac{1}{2}\sum_{i=1}^{m} \mathbf{R}_i^2 = \tfrac{1}{2}\|\mathbf{R}\|^2 \tag{3.17}$$

A number of numerical methods exist to solve the non-linear optimization problem (3.17), including steepest decent method, Newton method, Gauss-Newton method, Newton method, and Levenberg-Marquardt method[52,53]. In the current work, the Levenberg-Marquardt algorithm[54] and an adaptive nonlinear least-squares algorithm[55] are used for solving the unconstrained minimization problem (3.17). Interested readers can refer to these references on this topic for more details.

### 3.5 Mesh management strategies

As can be seen from the above two lemmas, the MDG+ICE method has the ability to resolve and fit discontinuities exactly by satisfying the interface condition. This is achieved by moving the grids, which are treated as an independent variable in the MDG+ICE formulation, such that their faces are aligned with solution discontinuities. Unfortunately, moving meshes inevitably become highly distorted, stretched, and tangled. Although the MDG+ICE method can handle highly distorted and even tangled grids where elements with zero or even negative volume exist, solutions obtained are either non-physical or meaningless. For example, a smooth multiple-valued solution, instead of a discontinuous solution, to the inviscid Burger's equation is numerically obtained by our MDG+ICE method, even though the Jacobian-determinant of some elements becomes negative in this case[49]. What is observed in the MDG+ICE solution process is that as the nonlinear residual decreases, elements in the vicinity of unfit interfaces are driven to flatten and eventually collapse. Elements are considered degenerate when the determinant of their Jacobian becomes negative or falls below a certain tolerance. By eliminating such degenerate elements, the discrete grid topology is modified to accommodate the interface topology present in the exact solution and ultimately track and fit these interfaces exactly. In general, all current existing strategies extensively developed for optimizing meshes can be used here to remove these degenerate elements. Degenerate triangles have circumcircles larger than their shortest edges, which can be classified as either a needle, whose longest edge is much longer than its shortest edge, or a cap, which has an angle close to 180° as shown in Figure 6a. Note that the classifications are not mutually exclusive, as shown in Figure 6b.

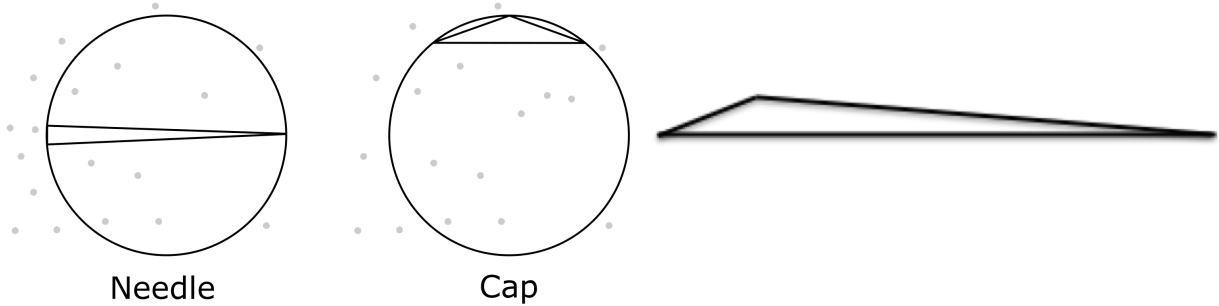

**Figure 6.** (a) (left) Classification of degenerate triangles and (b) (right) a triangle can be considered as either a cap or a needle.

Three popular grid operations: edge swap, edge collapse, and edge split as shown in Figure 7(a), are implemented in our MDG+ICE method in order to remove degenerate triangles, improve the quality of moving meshes, and maintain the same number of triangles, points, and boundary faces. An edge swap operation is widely used in Delaunay triangulation and minmax triangulation, which is both simple and local. Moreover, it does not change the number of grid points and triangles in a mesh, a very desired property. Edge swap is adopted here to both improve grid quality and remove cap triangles, as illustrated in Figure 7(b), where the cap triangle in green is removed as a result of an edge swap operation. Edge swap alone is mainly designed to improve the quality of triangles, and should not be allowed across interfaces. However, edge swap with a cap triangle removal is used everywhere and actually very desired and even necessary for capturing and fitting discontinuities. An edge collapse operation is very effective for removing needle triangles, although it cannot be directly used to eliminate cap triangles. Note that using meshing slicing[56], cap triangles can be removed one by one by splitting them into needles, which is not implemented in the current work. One drawback of edge collapse operations is the need to modify grids. In our implementation, each time a boundary or internal edge is collapsed, an edge split operation on the longest boundary or internal edge is followed in order to maintain the same number of triangles, grid points, and boundary faces. Since the MDG+ICE is based on the space-time formulation, the values of conservative variables can be easily interpolated or transferred to the mesh of the new topology in a simple and straightforward way without any concern to the conservative interpolation in contrast to the standard finite volume and DG formulations, which need to address this issue carefully.

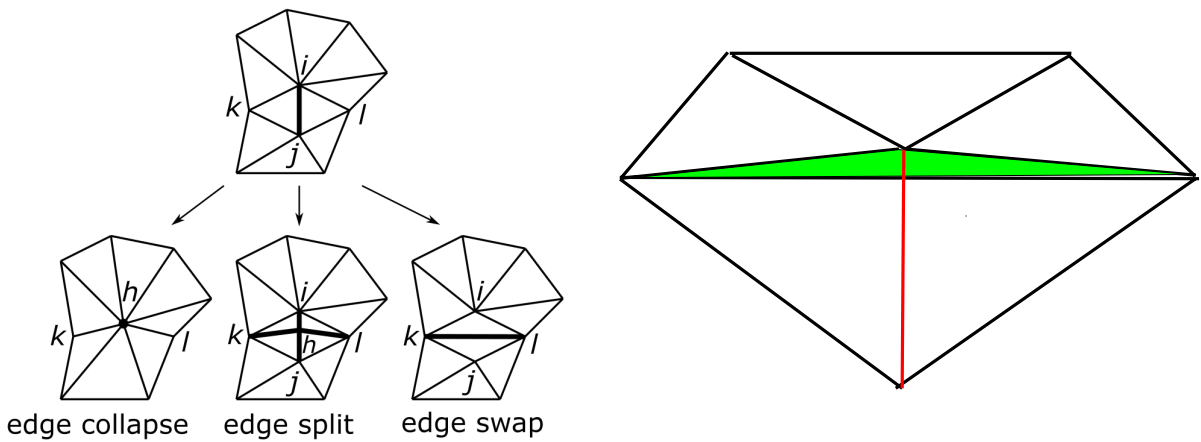

**Figure 7.** (a)(left) Illustration of three local grid operations and (b)(right) edge swap with a cap triangle removal

## IV. Numerical Examples

The developed MDG+ICE method is used to solve a variety of compressible flow problems. A few examples are presented here to demonstrate the accuracy, robustness, and ability of the MDG+ICE method for both 1D unsteady and 2D steady compressible Euler equations. Since there is no confusion, we will use abbreviation MDG+ICE and MDG interchangeably.

**A. Sod shock tube problem**
The classic Sod shock tube problem is one of the most widely used test cases, since it represents an exact solution to the full system of one-dimensional unsteady Euler equations. The exact solution contains simultaneously a shock wave, a contact discontinuity, and a rarefaction wave, which all emanate from one singularity point at t=0. This constitutes a particularly interesting and challenging problem for the space-time formulation where the singularity point exists and for the discontinuity tracking methods where different types of unknown interfaces exist. This example is chosen to assess (1) the ability of the MDG+ICE method to resolve singular points, and detect and fit discontinuities, (2) the robustness of the MDG+ICE method on highly distorted grids and grids with negative volumes, and (3) the order of $p$-convergence of the MDG+ICE method for discontinuous solutions. The initial conditions in the computation are the following:

$$(\rho, v, p)(x, t = 0) = \begin{cases} (1, 0, 1), & -0.5 \leq x < 0, \\ (0.125, 0, 0.1) & 0 < x \leq 0.5 \end{cases}$$

Computation is performed on a single space-time slab $\Omega(x, t) = (-0.5, 0.5) \times (0, 0.2)$ with 8 quadrilateral elements where the left 7 elements are initiated with the piecewise left constant state and the right one element is initiated with the piecewise constant right state. Numerical experiments are conducted on two grids, which are shown in Figure 8 along with the initial density fields. Six middle cells are degenerated at the origin (0,0), which are necessary to resolve the singularity point accurately and adequately. The second mesh contains an element with negative volume, as clearly indicated by the position of the initial interface. No significant difference between the MDG solutions obtained on these two grids is observed. Tangled meshes do not lead to a breakdown of the MDG solution process, indicating the robustness of the MDG method. Numerical solutions are computed using the MDG(P1), MDG(P2), and MDG(P3) methods. The results obtained by the MDG(P1), MDG(P2), and MDG(P3) solutions are presented in Figures 9, 10, and 11, respectively. Figures 9b, 10b, and 11b show a comparison of the density, pressure, velocity, and entropy production profiles at t=0.2 between the exact solution and the three computed MDG solutions, respectively. The MDG (P1) solutions are represented using a two-point connecting straight line on each element, while the MDG(P2) and MDG(P3) solutions are represented using a three-point connecting line on each element. As expected, all three MDG solutions are able to capture and fit both shock and contact as true discontinuities at the correct locations. In addition, the position of the head and tail of the refraction wave, where a derivative discontinuity exists, is also fit very accurately by the MDG(P1) solution and virtually resolved exactly by the MDG(P2) and MDG(P3) solutions. This superior accuracy of the MDG method should be attributed to the fact that the interface conservation is solved and enforced explicitly, as all discontinuity tracking and fitting methods are not purposely designed to accurately fit continuous solutions where their derivatives are discontinuous. A *p*-refinement study for this problem is conducted to numerically obtain quantitative measurement of the absolute errors and order of *p*-convergence for the MDG method. The following $L_2$ function norm in the space-time domain $\Omega(x,t)$ is used to measure the error of the MDG method

$$\left\| u_c - u_e \right\|_{L_2} = \sqrt{\int_{\Omega(x,t)} (u_c - u_e)^2 \, d\Omega}$$

where $u_c$ and $u_e$ are computed and exact solutions, respectively. The detailed results for the *p*-convergence study are presented in Table 1, where the number of degrees of freedom, the L2-error of the MDG solutions, and the order of convergence are listed. $L_2$ function norm of the error function is plotted in Figure 12 against the inverse of the square root of the number of degrees of freedom on a log-log scale. One can clearly observe that the MDG method exhibits an exponential rate of convergence for *p*-type refinement, which is reflected by the downward curving line instead of a straight line for the *h*-refinement. This example demonstrates some of the most attractive features of the MDG method due to the interface conservation enforcement: being able to detect and fit different types of discontinuities automatically without resorting to specialized logics specific to each type and achieving the formal rate of *p*-convergence even for discontinuous solutions and solutions with singularities, that no other numerical methods can to the best of our knowledge. Solving and enforcing the conservation laws on interfaces enable the MDG+ICE method to not only detect and fit the contact discontinuity and shock wave but also exactly resolve the head and tail of the rarefaction wave, where the solutions are continuous but their derivatives are discontinuous.

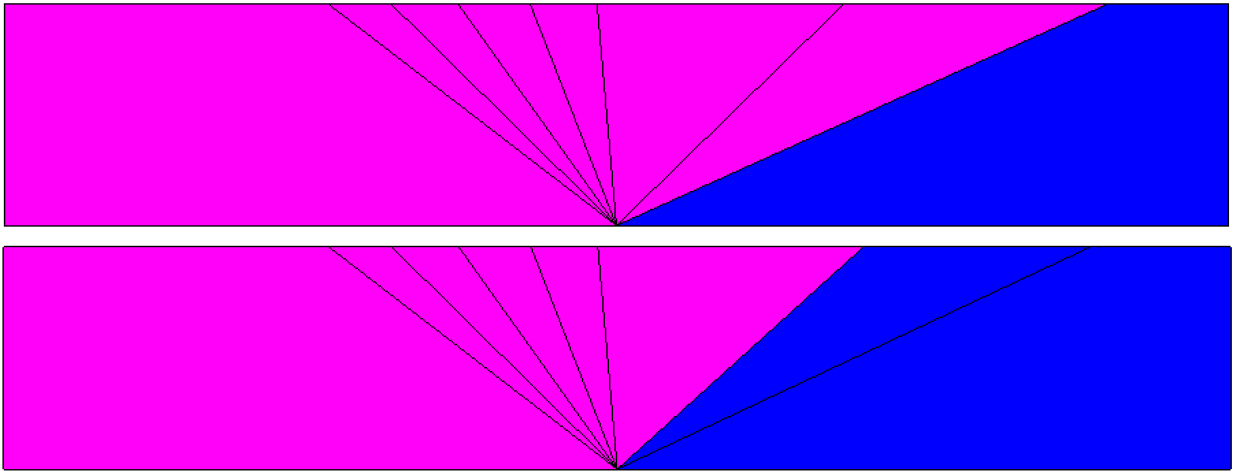

**Figure 8. Initial space-time grids and density fields on Ω=(-0.5,0.5)×(0,0.2) for Sod shock tube problem**

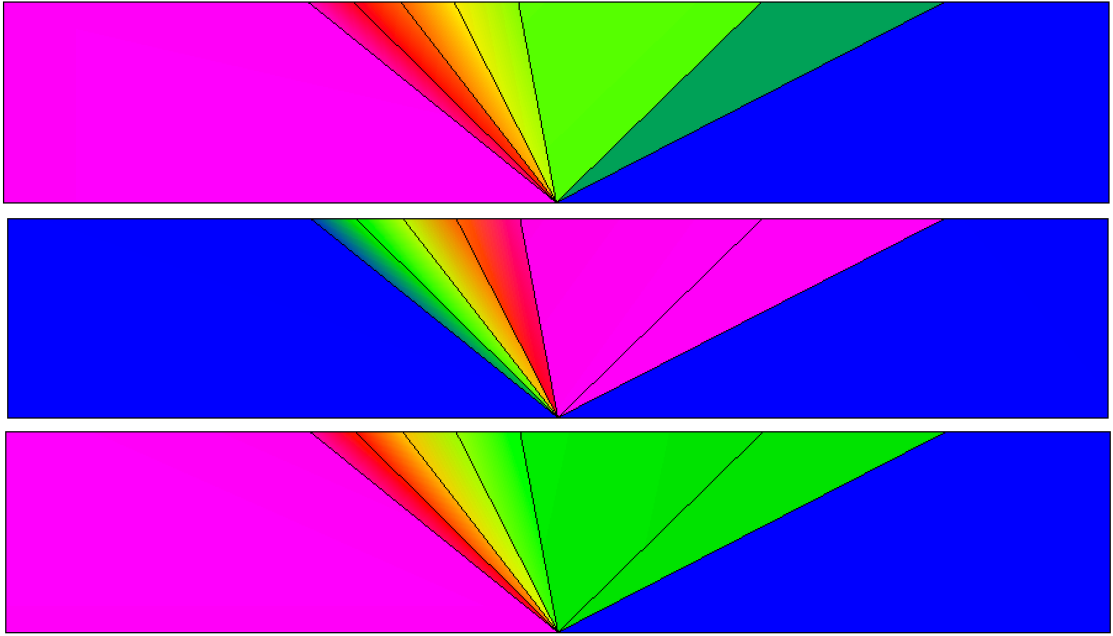

**Figure 9(a). Computed density (top), velocity (middle), and pressure (bottom) along with the converged spacetime grid obtained by the MDG(P1) method on Ω=(-0.5,0.5)×(0,0.2) for Sod shock tube problem**

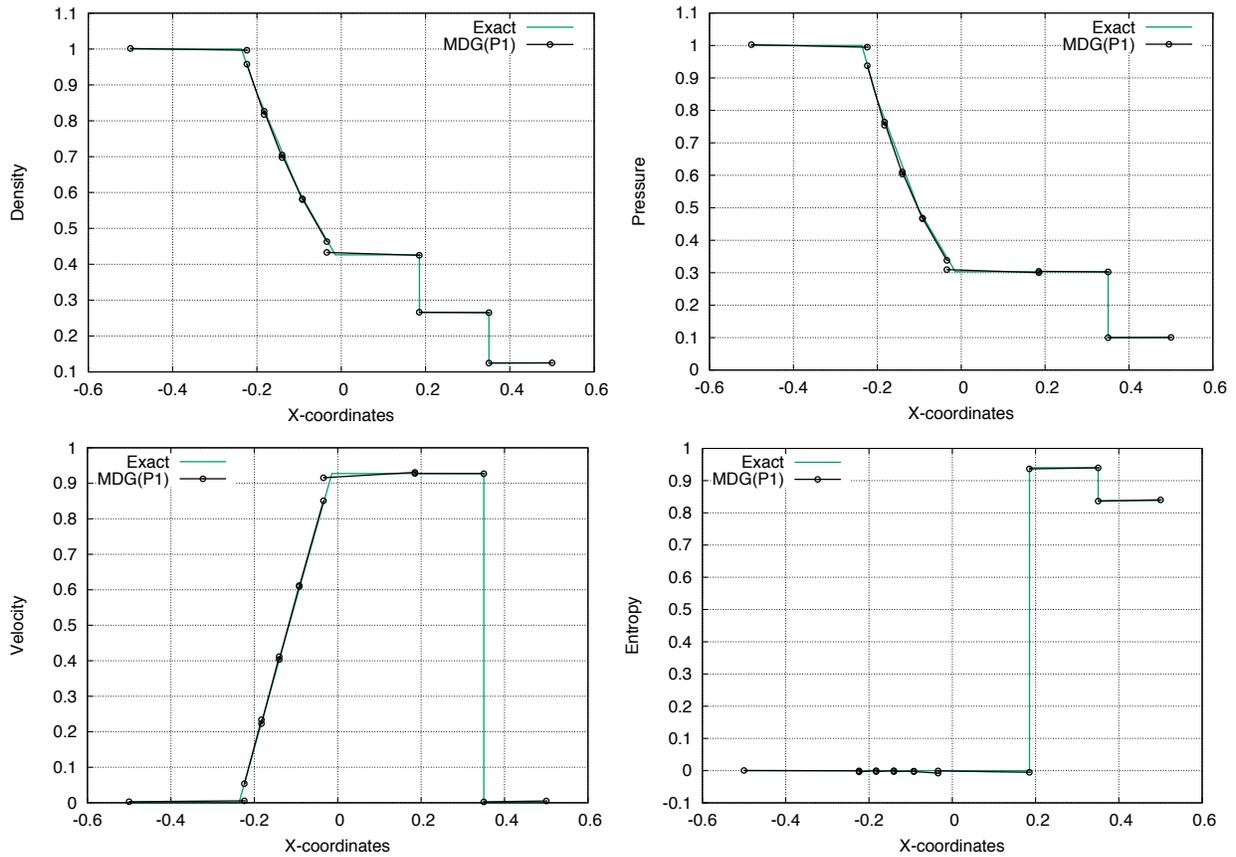

**Figure 9(b). Comparison of density, velocity, pressure, and entropy production at t=0.2 between the MDG(P1) solution and the analytical solution.**

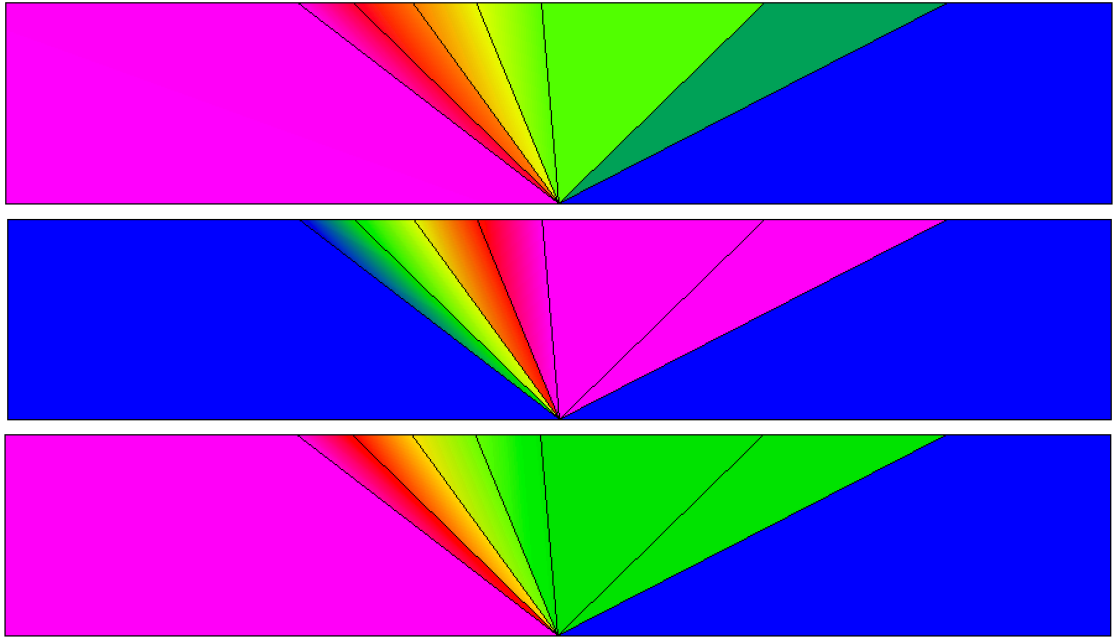

**Figure 10(a). Computed density (top), velocity (middle), and pressure (bottom) along with the converged space-time grid obtained by the MDG(P2) method on Ω=(-0.5,0.5)×(0,0.2) for Sod shock tube problem**

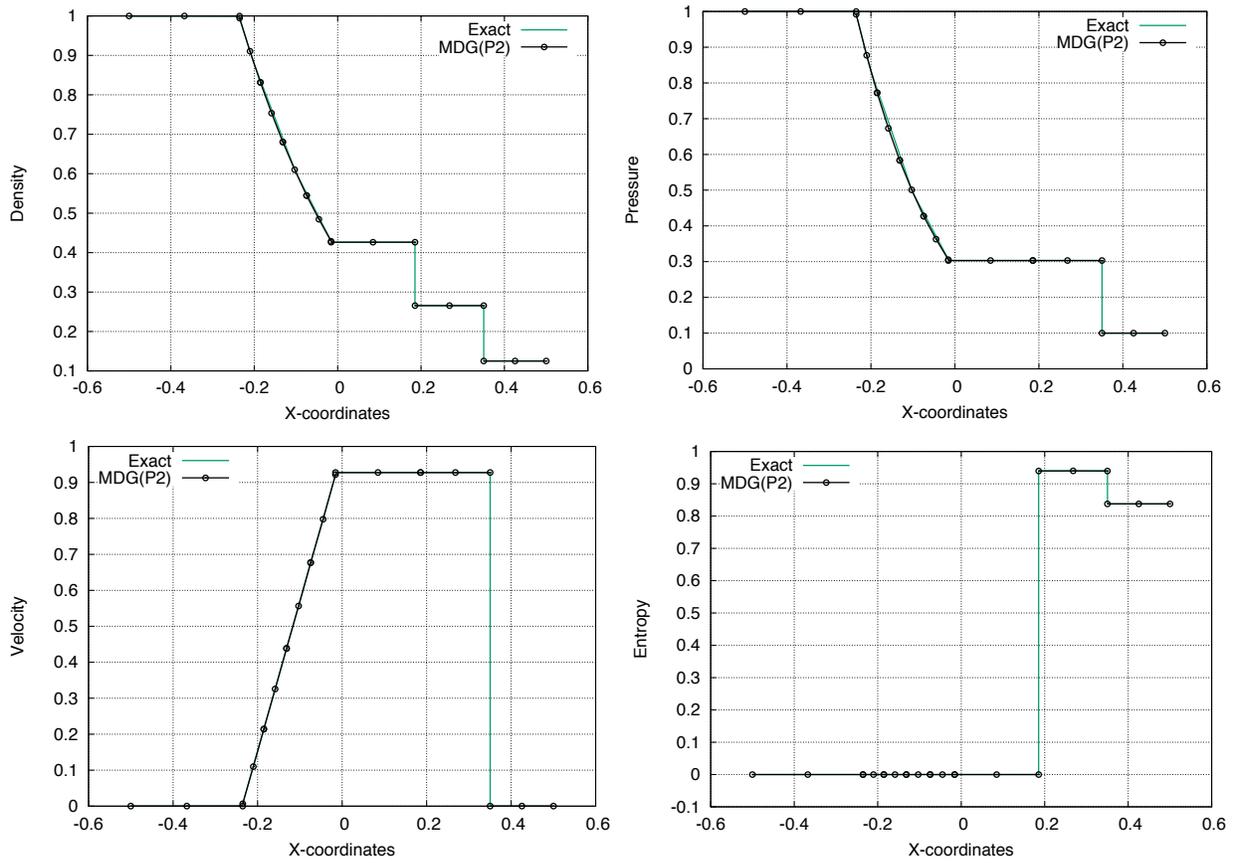

**Figure 10(b). Comparison of density, velocity, pressure, and entropy production at t=0.2 between the MDG(P2) solution and the analytical solution.**

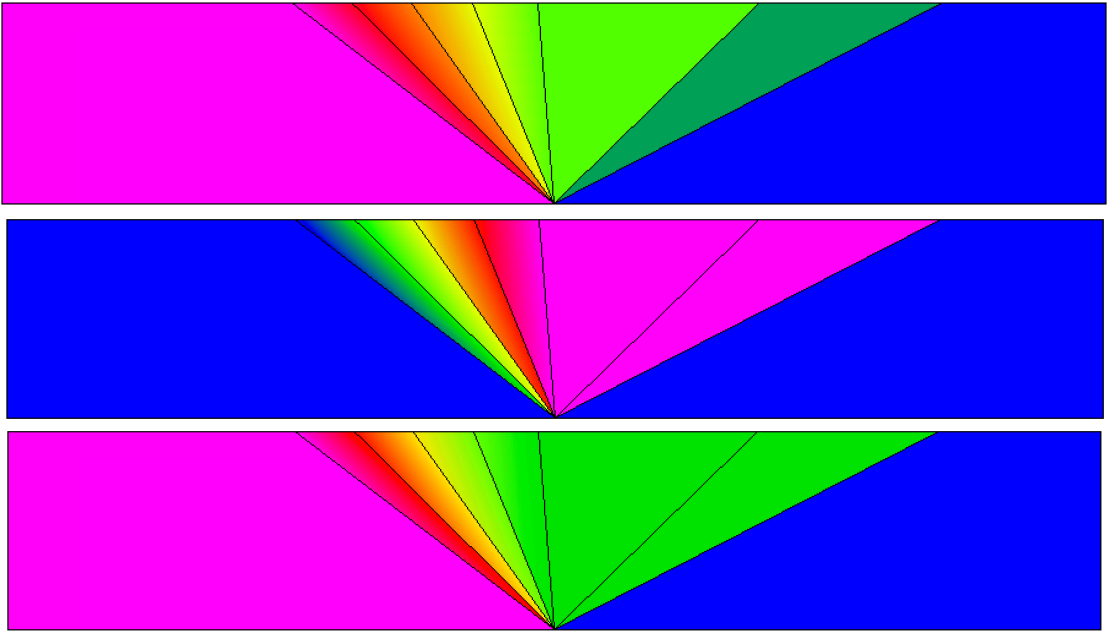

**Figure 11(a). Computed density (top), velocity (middle), and pressure (bottom) along with the converged spacetime grid obtained by the MDG(P3) method on Ω=(-0.5,0.5)×(0,0.2) for Sod shock tube problem**

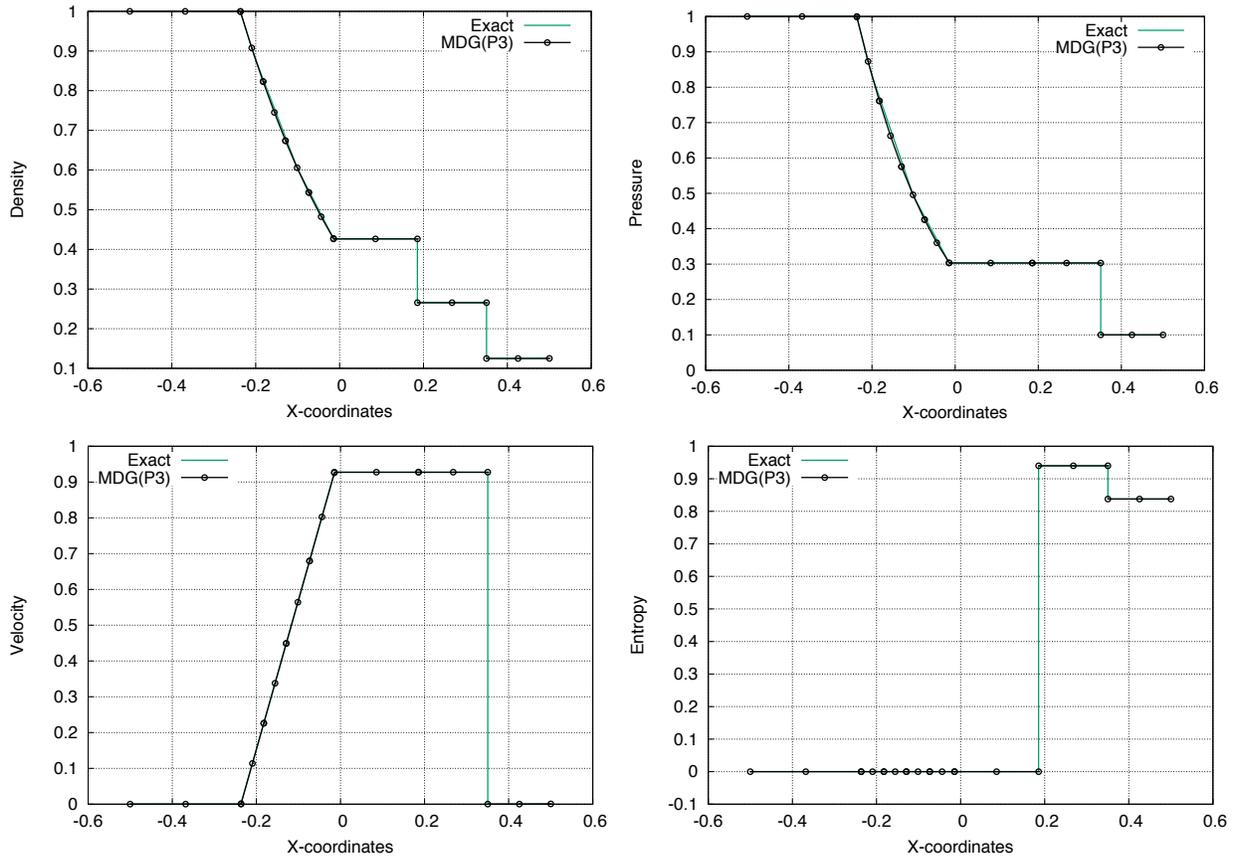

**Figure 11(b). Comparison of density, velocity, pressure, and entropy production at t=0.2 between the MDG(P3) solution and the analytical solution.**

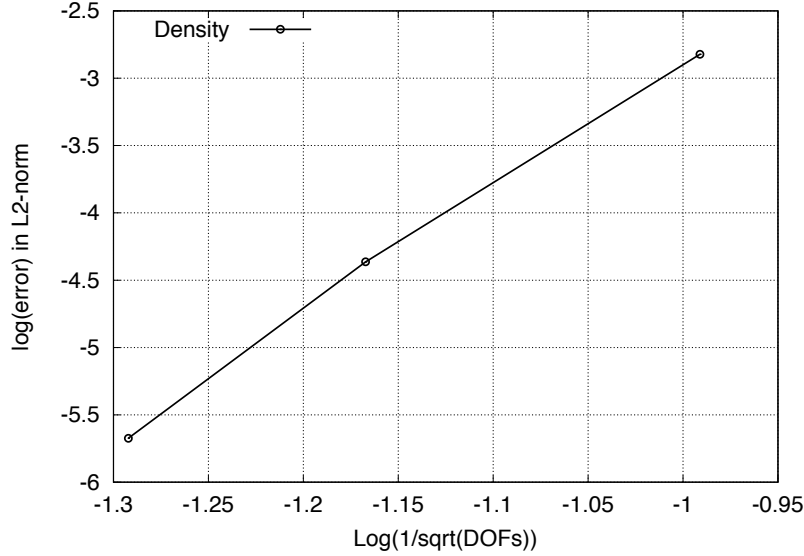
**Figure 12. Convergence histories of p-type refinement for the MDG methods**

**Table 1. L$_2$-error and order of convergence for the MDG methods.**

| MDG | Log(1/√DOFs) | Log($\|\rho - \rho_e\|_{L_2(\Omega)}$) | Slope |
|---|---|---|---|
| P1 | -0.9911 | -0.282255E+01 | |
| P2 | -1.1672 | -0.436399E+01 | 8.75 |
| P3 | -1.29216 | -0.567351E+01 | 10.48 |

**B. Lax–Harden Riemann problem**

This is another well-known test case for the shock tube problem. The initial conditions in the present computation are the following:

$$(\rho, v, p)(x, t = 0) = \begin{cases} (0.445, 0.698876404, 3.52773), & -0.5 \leq x < 0, \\ (0.5, 0, 0.571) & 0 < x \leq 0.5 \end{cases}$$

Computation is performed using 8 quadrilateral cells on a space-time domain $\Omega = (-0.5, 0.5) \times (0, 0.15)$ where the left 7 elements are initiated with the piecewise left constant state and the right one element is initiated with the piecewise constant right state. Initial mesh and density field are shown in Figure 13. Six middle cells are degenerated at the origin (0,0), which are necessary to resolve the singularity point accurately and adequately. Numerical solutions are computed using the MDG(P1), MDG(P2), and MDG(P3) methods. The results obtained by the MDG(P1), MDG(P2), and MDG(P3) solutions are presented in Figures 14, 15, and 16, respectively. Figures 14b, 15b, and 16b show a comparison of the density, pressure, velocity, and entropy production profiles at t=0.2 between the exact solution and the three computed MDG solutions, respectively. The converged space-time grids are different for the three MDG methods, due the fact that the location of grid points is not uniquely determined in the constant flow region. However, all three MDG solutions are able to fit both shock and contact as true discontinuities at the correct locations. Moreover, the position of the head and tail of the refraction wave, where a derivative discontinuity exists, is also fit accurately by the MDG(P1) solution due to lack of resolution and virtually exactly by the MDG(P2) and MDG(P3) solutions. As in the previous example, a *p*-refinement study for this problem is conducted to numerically obtain quantitative measurement of the absolute errors and order of *p*-convergence for the MDG method. The detailed results for the *p*-convergence study are presented in Table 2, where the number of degrees of freedom, the L2-error of the MDG solutions, and the order of convergence are listed. L$_2$ function norm of the error function is plotted in Figure 17 against the inverse of the square root of the number of degrees of freedom on log-log scale. An exponential rate of convergence for *p*-type refinement is observed again in this test case, confirming again that the MDG methods are able to detect and track multiple types of discontinuities automatically and achieve the formal rate of *p*-convergence even for discontinuous solutions and solutions with singularities.

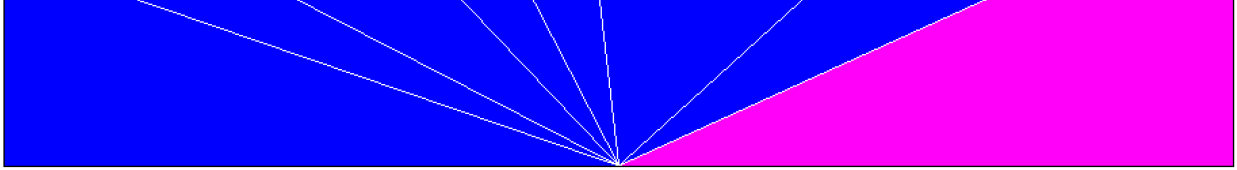

**Figure 13. Initial space-time grid and density field for Lax-Harden Riemann problem**

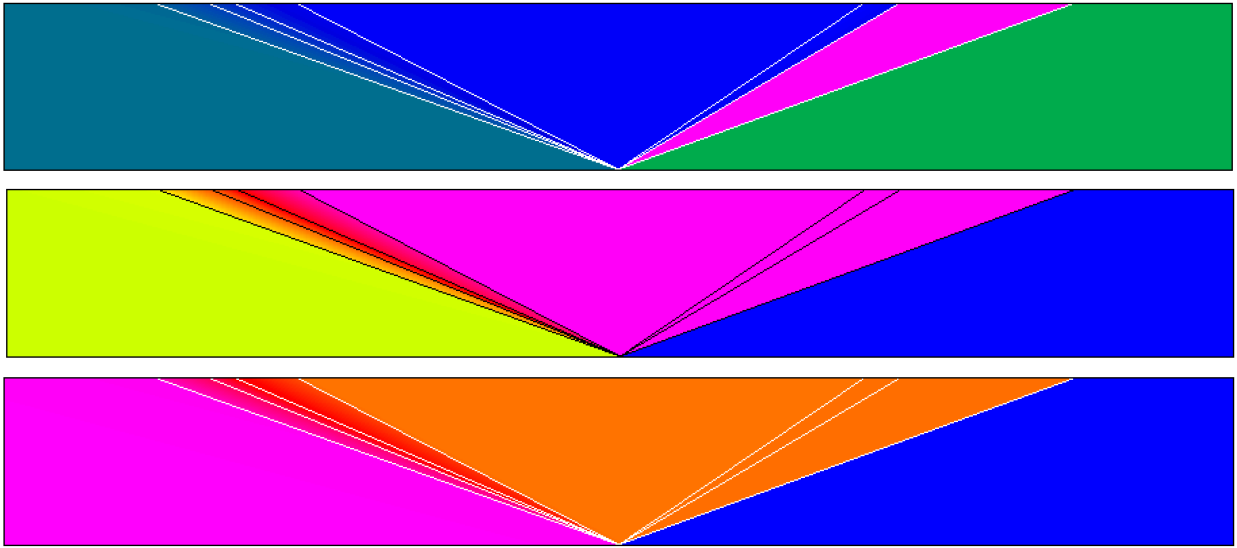

**Figure 14(a). Computed density (top), velocity (middle), and pressure (bottom) along with the converged spacetime grid obtained by the MDG(P1) method on Ω=(-0.5,0.5)×(0,0.15) for Lax-Harden Riemann problem**

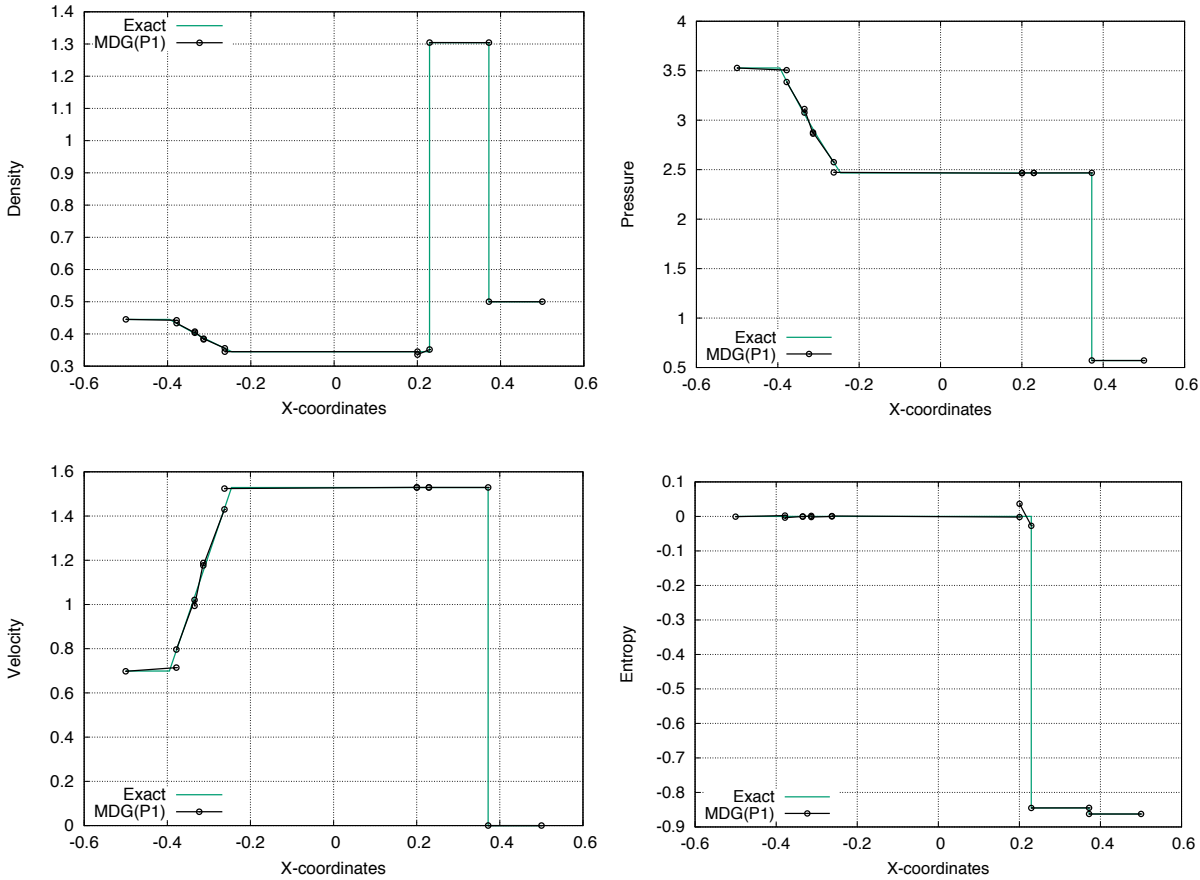

**Figure 14(b). Comparison of density, velocity, pressure, and entropy production at t=0.15 between the MDG(P1) solution and the analytical solution.**

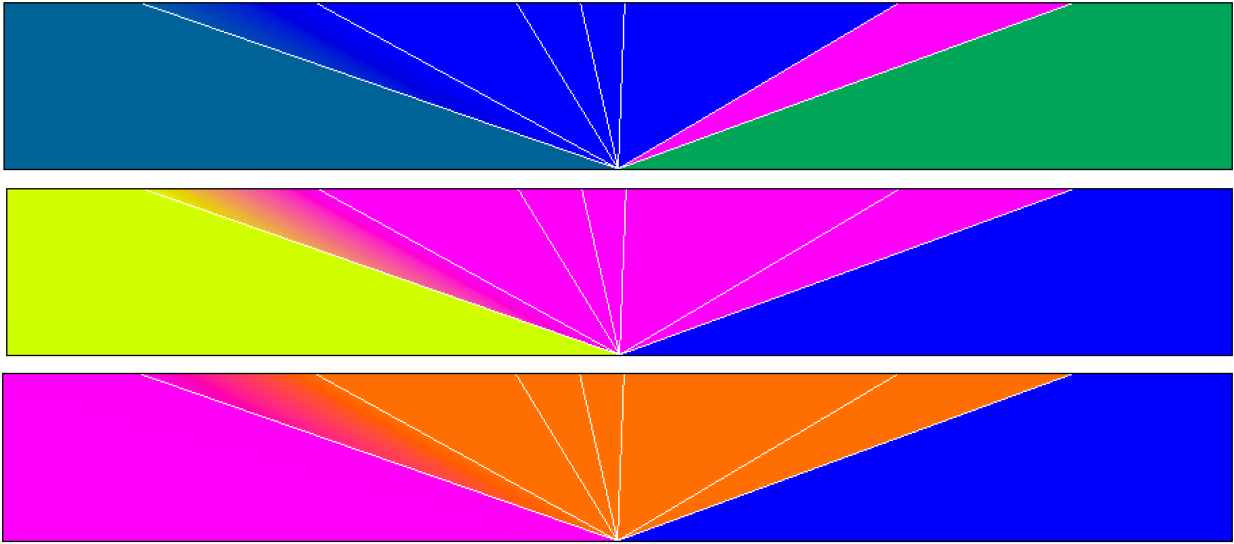

**Figure 15(a). Computed density (top), velocity (middle), and pressure (bottom) along with the converged spacetime grid obtained by the MDG(P2) method on Ω=(-0.5,0.5)×(0,0.15) for Lax-Harden Riemann problem**

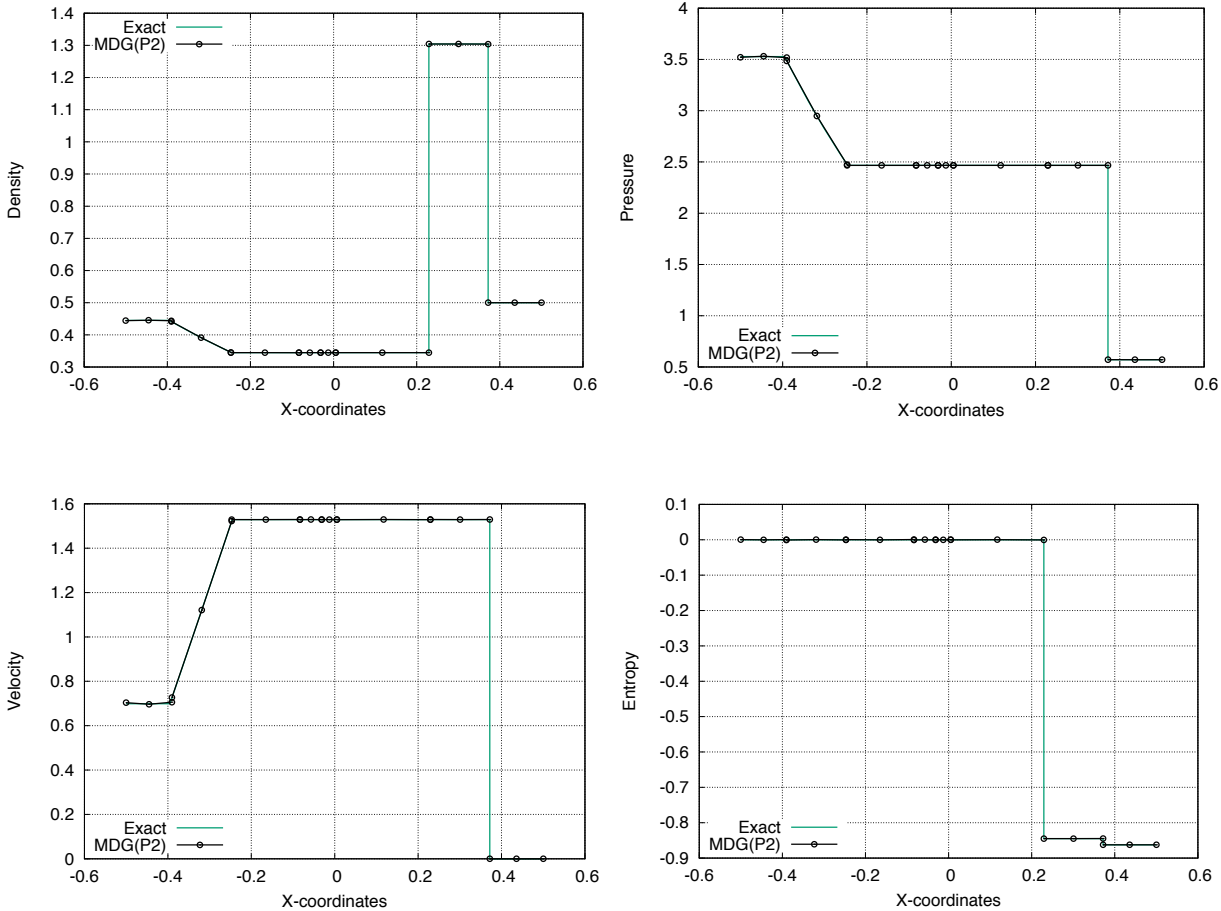

**Figure 15(b). Comparison of density, velocity, pressure, and entropy production at t=0.15 between the MDG(P2) solution and the analytical solution.**

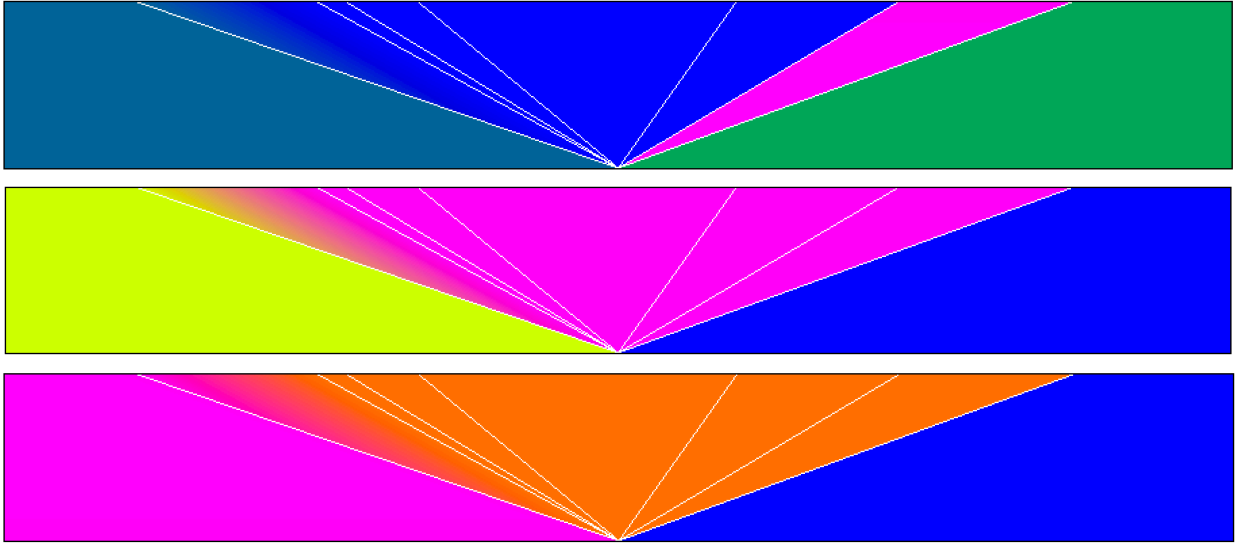

**Figure 16(a). Computed density (top), velocity (middle), and pressure (bottom) along with the converged spacetime grid obtained by the MDG(P3) method on Ω=(-0.5,0.5)×(0,0.15) for Lax-Harden Riemann problem**

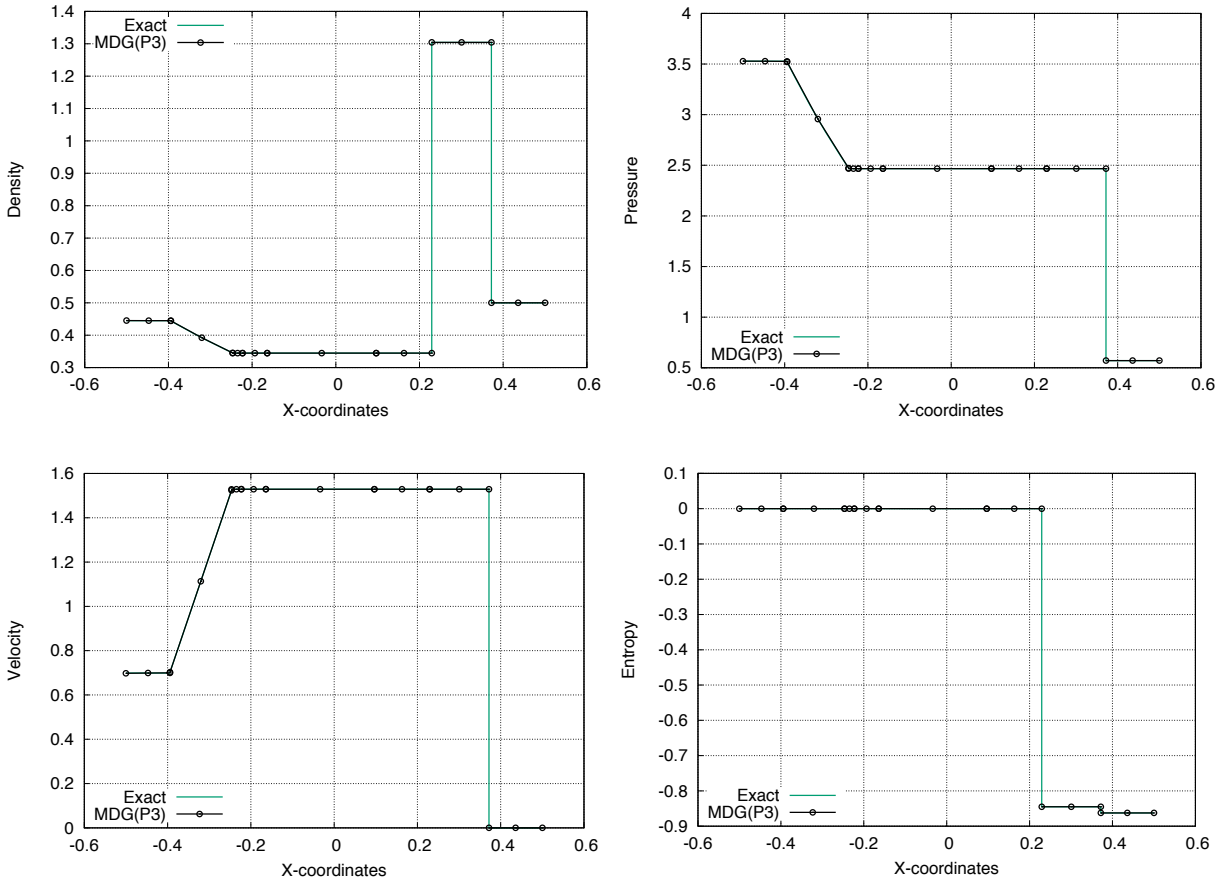

**Figure16(b). Comparison of density, velocity, pressure, and entropy production at t=0.15 between the MDG(P3) solution and the analytical solution.**

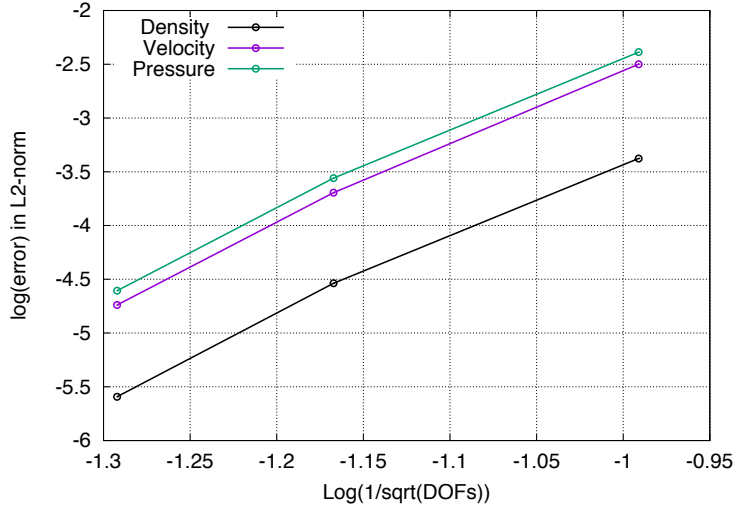

**Figure 17. Convergence histories of density, velocity, and pressure for p-refinement**

**Table 1. L$_2$-error and order of convergence for the MDG methods.**

| MDG | Log($1/\sqrt{\text{DOFs}}$) | Log($\|\rho - \rho_e\|_{L_2(\Omega)}$) | Slope |
|---|---|---|---|
| P1 | -0.9911 | -0.337635E+01 |  |
| P2 | -1.1672 | -0.453702E+01 | 6.59 |
| P3 | -1.29216 | -0.559210E+01 | 8.44 |

**C. Receding flow problem - 123 problem**

The receding flow problem, also known as the 123 problem or the double rarefaction wave problem, consists of a smooth flow undergoing rarefaction caused by two flows receding from each other. The receding problem was extensively studied and identified by Liou[57] as an open numerical problem, as all well-known numerical methods in different formulations (Eulerian, Lagrangian, and ALE) produce an anomalous temperature rise, often termed "overheating", at the origin that cannot be fixed by using different numerical fluxes, refining the mesh, decreasing the time-step, or increasing the order of numerical methods. In the receding flow, the initial discontinuous solution (t=0) singularly becomes a continuous solution when t>0, which causes the overheating. In fact, all numerical methods can trivially solve this problem, if started from the exact solution at t=$\varepsilon > 0$. However, none of them can survive in the time window [0,$\varepsilon$], as Liou showed that the overheating occurs in the very first instants when the two receding rarefaction waves introduce the most significant discontinuities. The setup of the problem in this numerical experiment is the following:

$$(\rho, v, p)(x, t = 0) = \begin{cases} (1, -2, 0.4), & \forall -0.5 \leq x < 0, \\ (1, 2, 0.4) & \forall 0 < x \leq 0.5 \end{cases}$$

A mesh of 16 quadrilateral cells on a space-time domain $\Omega$ =(-0.5,0.5)×(0,0.15) shown in Figure 18 is used to solve this challenging problem, where the left 8 elements are initiated with the piecewise left constant state and the right 8 elements are initiated with the piecewise constant right state. 14 middle cells are degenerated at the origin (0,0), which are necessary to resolve the singularity point accurately and adequately. Numerical solutions are computed using the MDG(P2) method. This test case is chosen to investigate if the MDG method is able to solve this open numerical problem and demonstrate the robustness of the MDG+ICE method in the sense that elements with zero volume and even negative pressure and density do not lead to a breakdown of the solution process in the MDG+ICE method. Contrary to Sod and Lax-Harden problems where discontinuities need to be captured and fitted, the challenge for the receding problem is the need to remove the initial discontinuity and accurately resolve the smooth solution. In the MDG method, this is automatically achieved by the collapse of the two quadrilateral cells adjacent to the initial interface, as can be observed in Figure19, where the converged space-time grid and velocity field are presented. Figure 20 shows a comparison of the density profile at t=0.15 between the computed MDG(P2) solution and analytical solution. One can observe the two zero-volume elements and the negative values of the density in the two collapsed cells, which is intentionally left in the graph to demonstrate that elements with zero volumes and negative density do not break down the solution process. In practice, these cells with zero and negative-volumes are simply removed by the mesh management strategies presented in section 3.5. Figure 21 presents a comparison of the density, pressure, velocity, and internal energy profiles at t=0.15 between the exact solution and the computed MDG solution after the removal of the two collapsed cells. Astoundingly, a highly accurate solution without overheating is obtained by the MDG+ICE method, which no other numerical methods can to the best of our knowledge. This achievement is attributed to the ability of the MDG+ICE method to remove the singularity solution via the interface condition enforcement and the grid movement and management.

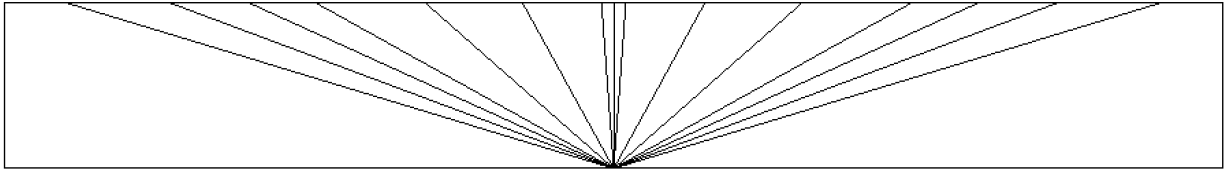

**Figure 18. Initial space-time grid on $\Omega$=(-0.5,0.5)×(0,0.15) for the receding flow problem**

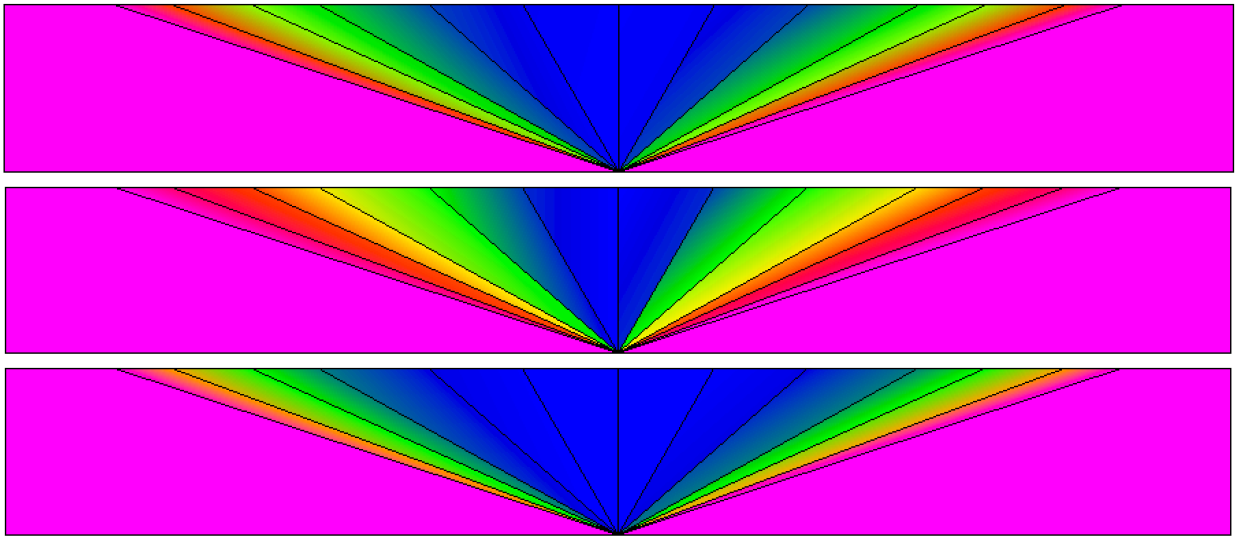

**Figure 19. Converged space-time grid and computed density, velocity, and pressure obtained by the MDG(P2) method for the receding flow problem**

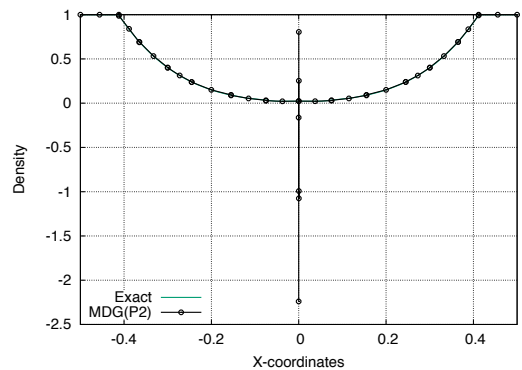

**Figure 20. Comparison of the density profile at t=0.15 between the MDG(P2) solution and the analytical solution**

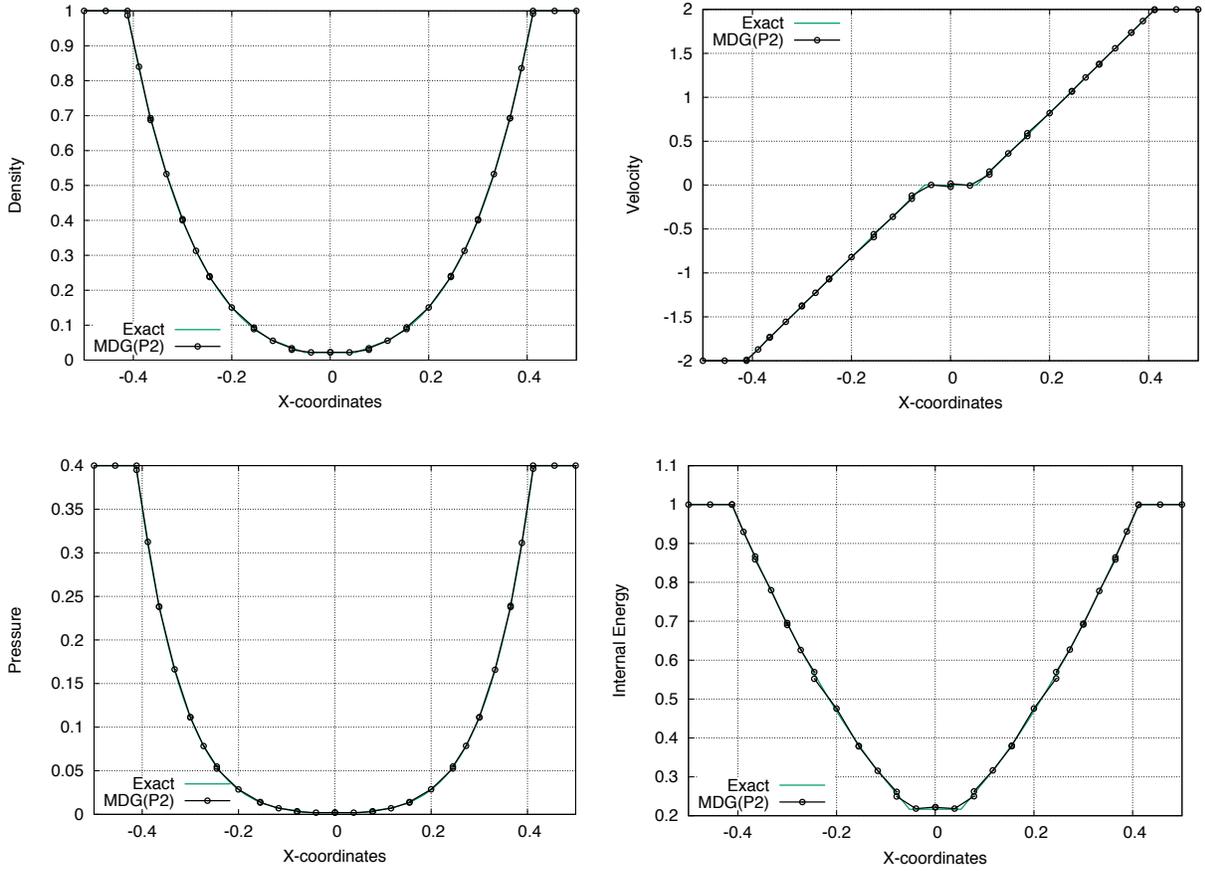

**Figure 21. Comparison of density, velocity, pressure, and internal energy at t=0.15 between the MDG(P2) solution and the analytical solution**

**D. Noh problem**

The 1D Noh problem is a well-studied and widely used test problem to assess the accuracy and robustness of the different numerical schemes. The solution of the Noh problem consists of two infinite strength shocks moving out from the center, leaving a constant density and pressure state behind. Similar to the receding problem, the Noh problem represents a great challenge to all numerical schemes in all (Eulerian, Lagrangian, and ALE) formulations, as they all unavoidably produce an anomalous temperature rise, the so-called overheating, at the origin. The setup of the problem is as follows:

$$(\rho, v, p)(x, t = 0) = \begin{cases} (1, 1, 1.0E - 06), & \forall - 0.5 \leq x < 0, \\ (1, -1, 1E - 06) & \forall 0 < x \leq 0.5 \end{cases}$$

The ratio of the specific heats is 5/3 instead of 1.4 in this test case. A mesh of 4 quadrilateral cells on a space-time domain $\Omega =(-0.5,0.5)\times(0,1)$ is used to solve this challenging problem, where the left two elements are initiated with the piecewise left constant state and the right 2 elements are initiated with the piecewise constant right state. 2 middle cells are degenerated at the origin (0,0), which are necessary to resolve the singularity point accurately and adequately. Numerical solutions are computed using the MDG(P1) method. This test case is chosen to investigate if the MDG method is able to solve a very strong shock wave problem and demonstrate the robustness of the MDG+ICE method in the sense that elements with zero volume and even negative pressure and density do not break down the solution process in the MDG+ICE method. Similar to the previous receding flow problem, one element is automatically collapsed and therefore removed during the process of the nonlinear iterations, as can be observed in Figure 22, where the initial space-time mesh and the computed space-time grid and density field at different non linear iterations are presented. One also observed that the density and pressure become negative during the course of the solution process, which nevertheless does not have any adverse effect on the accuracy of the final solution.

Figure 23 illustrates the convergence histories of the nonlinear residuals for the conservation laws and the interface conservation, where $L_2$ vector norm of the nonlinear residuals is plotted against the number of nonlinear iterations. A good convergence history is observed even for this challenging problem. Figure 24 presents a comparison of the density, pressure, velocity, and internal energy profiles at t=1 between the exact solution and the computed MDG solution. The computed MDG(P1) solution is virtually identical to the exact solution, and does not exhibits any anomalous overheating that exists in all other numerical methods, clearly demonstrating the robustness and accuracy of the MDG+ICE method for this challenging problem.

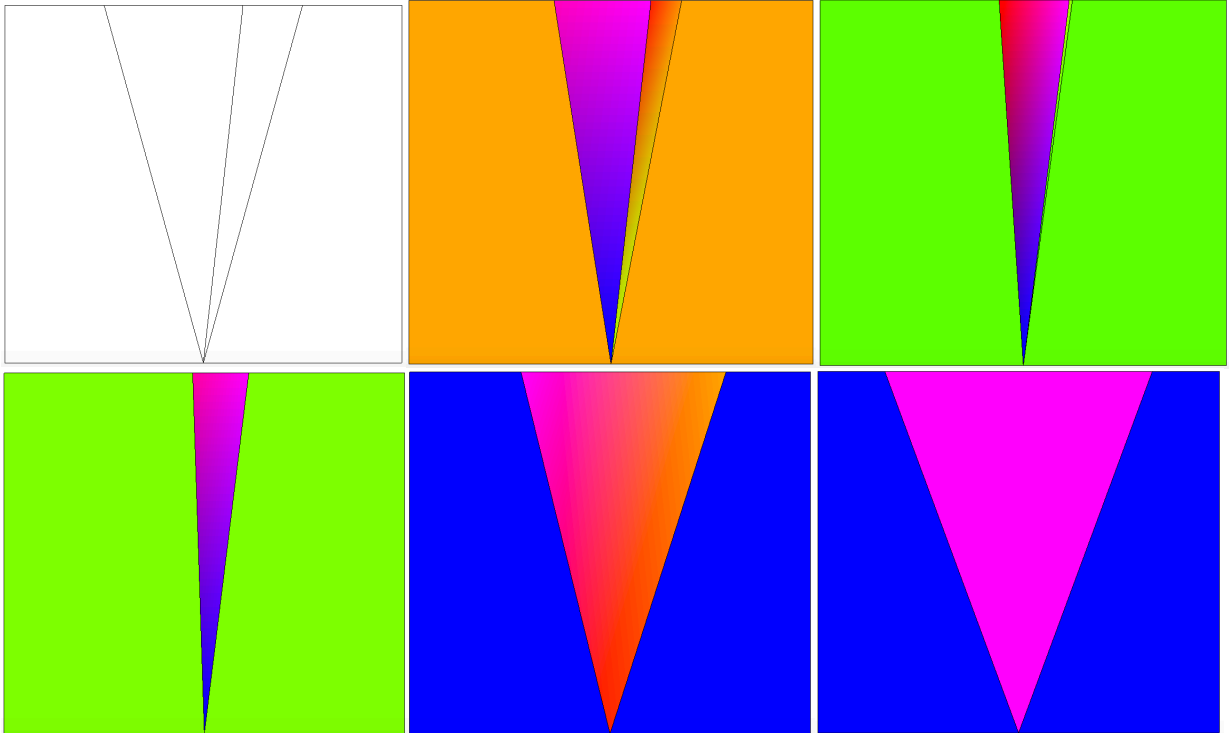

**Figure 22. Initial grid and computed space-time grids and density fields at different iterations (from left to right and top to bottom: 0, 2, 10, 20, 30, 40) by the MDG(P1) method for the Noh problem**

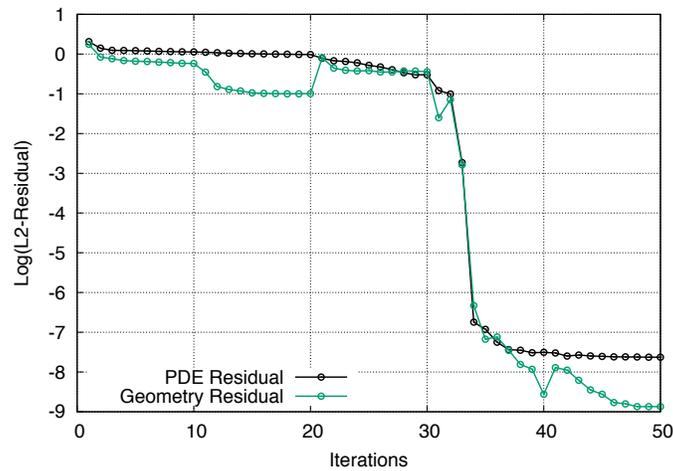

**Figure 23. Convergence histories of nonlinear residuals for conservation laws and interface conservation**

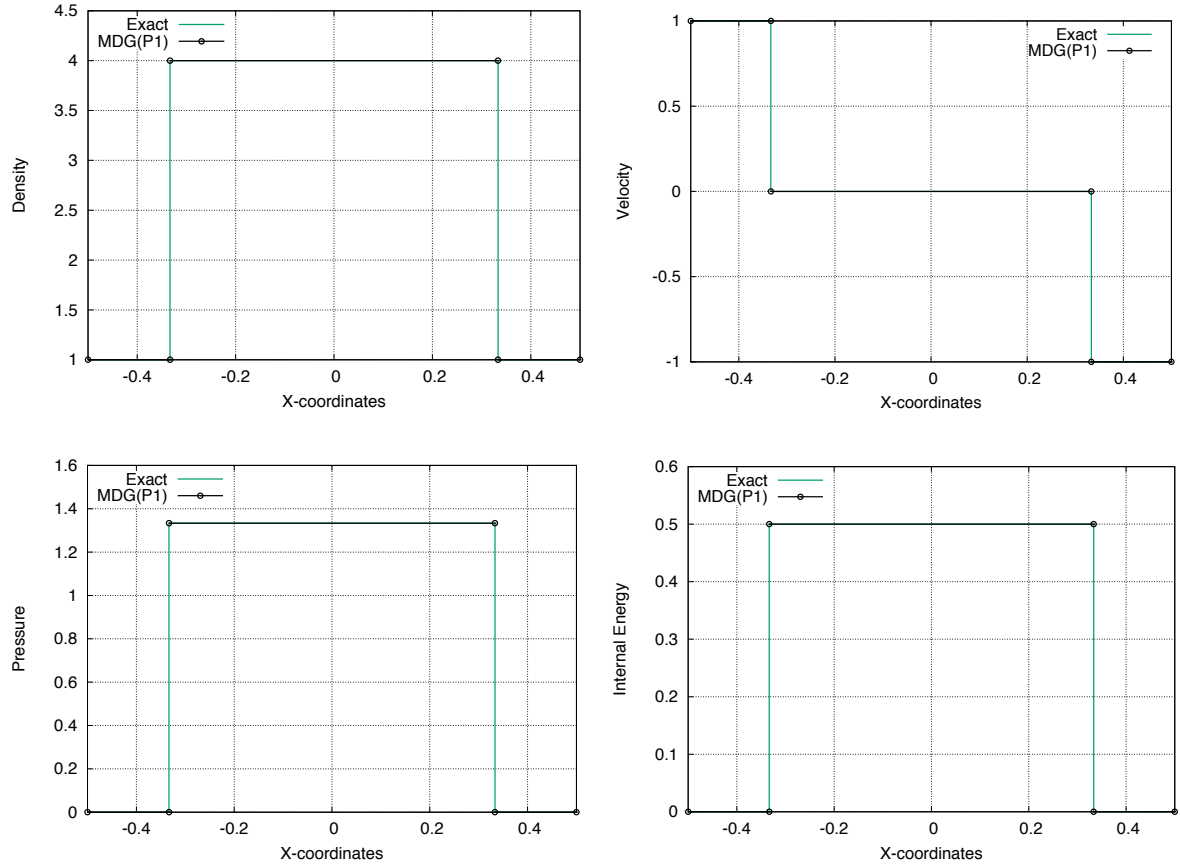

**Figure 24. Comparison of density, velocity, pressure, and internal energy at t=1 between the MDG(P1) solution and the analytical solution for the Noh problem.**

**D. Supersonic flow past a wedge.**

A supersonic flow at $M_\infty = 3$ past a 15° wedge is considered in this example. A numerical solution obtained by the standard DG(P0) method is used as the initial condition for the MDG(P1) method. The initial mesh has 189 elements, 115 points, and 39 boundary faces. The computed meshes and the density fields at different iterations are shown in Figure 25. As the nonlinear iterations progress, triangles in the vicinity of the shock wave are driven to flatten and collapse, and consequently meshes will be inevitably become highly distorted, stretched, and tangled. The highly distorted and degenerated elements are effectively removed from the grid by the mesh management strategies described in Section 3. One can clearly see that the oblique shock is captured and fitted as a true discontinuity with an instantaneous jump across the shock interface. Figure 26(a) presents the convergence histories of the nonlinear residuals for the conservation laws, the interface condition, and the combination of both, where $L_2$ vector norm of the nonlinear residuals is plotted against the number of nonlinear iterations. One can observe that the residuals for both conservation laws and interface condition converge at the same rate. The computed density distribution along a horizontal line at a distance of 0.3625 from the bottom is compared with the one from the analytical solution in Figure 26(b). The MDG(P1) solution and analytical solution are virtually identical. This example demonstrates the ability of MDG+ICE method to resolve, detect, and track an initially unfitted shock.

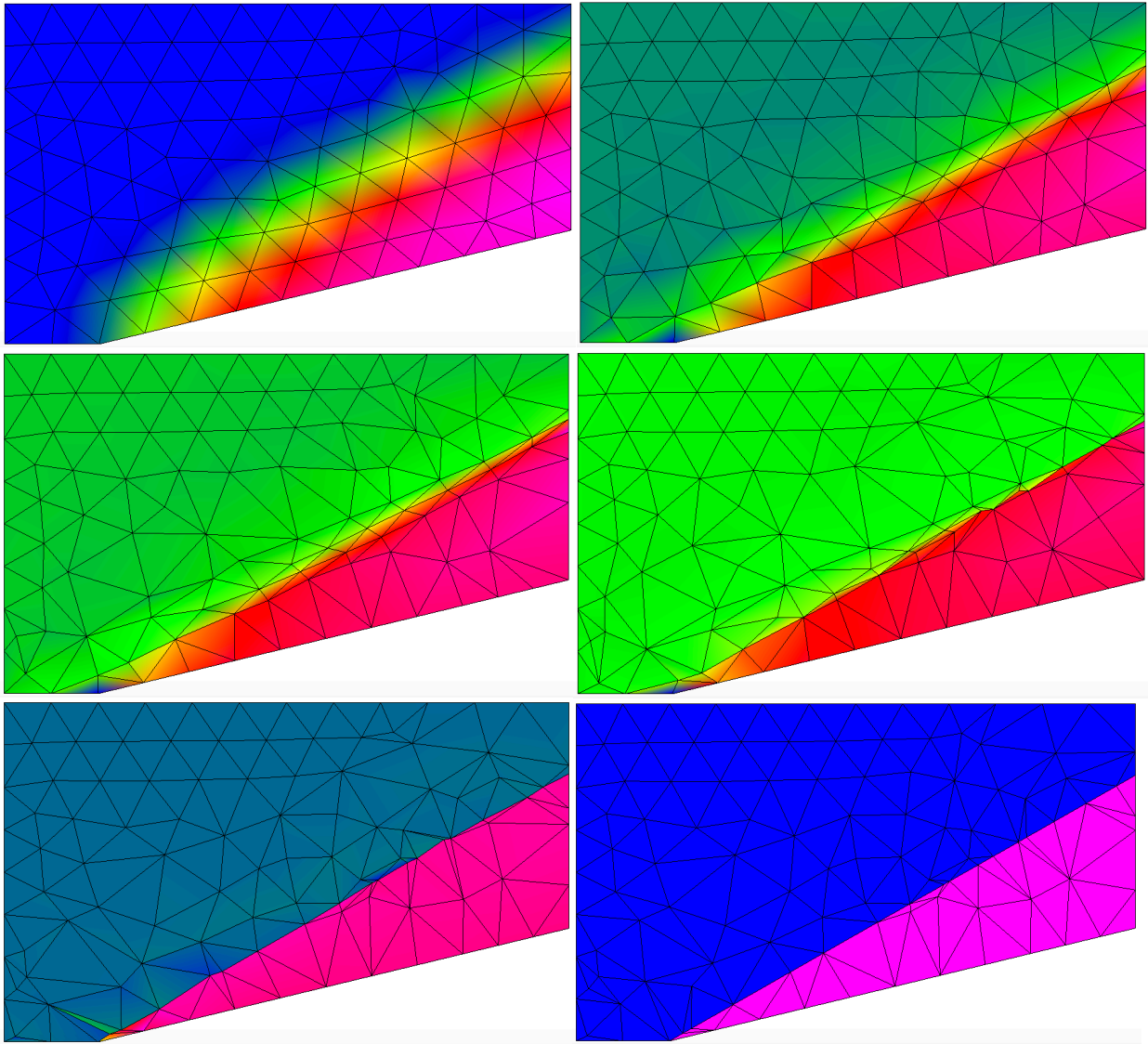

**Figure 25. Computed meshes and density fields at different iterations (from left to right and top to bottom: 0, 4, 8, 12, 16, 20) for a supersonic flow past a wedge.**

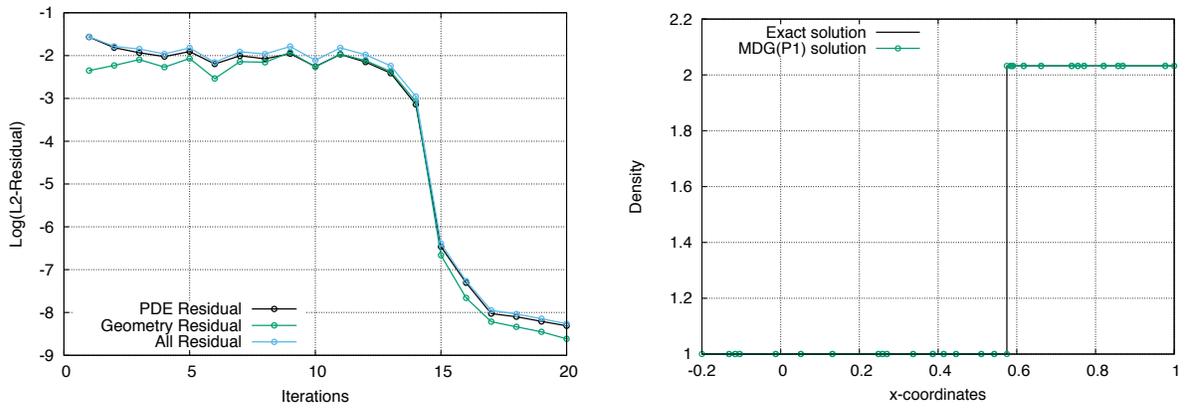

**Figure 26. (a) (left) Convergence histories of nonlinear residuals for conservation laws and interface conservation and (b)(right) comparison of the density profile at y=0.3625 between the MDG(P1) solution and exact solution.**

**D. Supersonic flow past a diamond wedge.**

A supersonic flow at $M_\infty = 2$ past a diamond wedge in a channel is considered in this example. This example is chosen to assess the ability of the MDG+ICE method for fitting discontinuities with non-trivial topology. The mesh consists of 500 elements, 284 points, and 66 boundary faces. The density fields obtained by the DG(P0) solution on the initial grid and the MDG(P1) solution on the final grid are shown in Figure 27. One can observe that shock reflection, shock-shock interaction, and expansion fan are well resolved by the MDG+ICE method, clearly demonstrating the ability of the MDG+ICE method to detect, resolve, and fit discontinuities with complex topology.

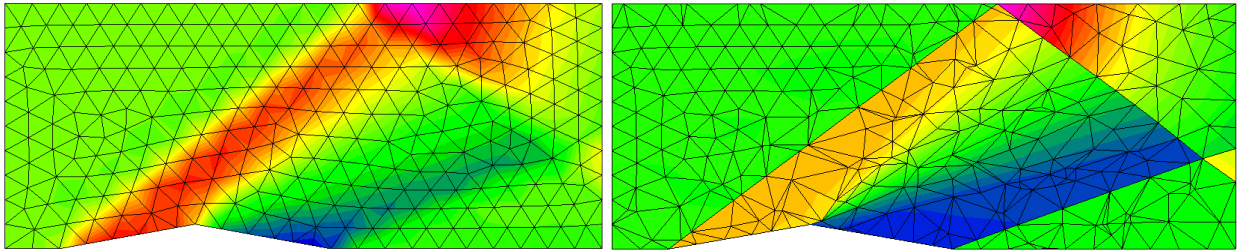

**Figure 27. Density field computed by the DG(P0) solution on the initial grid (left) and density fields obtained by the MDG(P1) solution on the final grid for a supersonic flow past a diamond wedge**.

**E. Hypersonic flow past a half-cylinder**

Two hypersonic flows past a half-cylinder at a free stream Mach number of $M_\infty = 8$ and 15, respectively, are considered in this example. This example is chosen to assess the ability of the MDG+ICE method to fit extremely strong shock waves. The mesh used for the Mach number 8 flow has 315 elements, 191 points, and 65 boundary faces, while the mesh used the hypersonic flow at a free-stream Mach number of 15 has 497 elements, 284 points, and 69 boundary faces. The initial mesh and density field obtained by the DG(P0) solution, and the converged mesh and density field obtained by the MDG(P1) solutions are shown in Figure 28 for these two test cases. One does not observe the appearance of any shock anomalies such as instabilities and carbuncle phenomena, clearly demonstrating the ability of the MDG+ICE method to detect, track, and fit strong discontinuities.

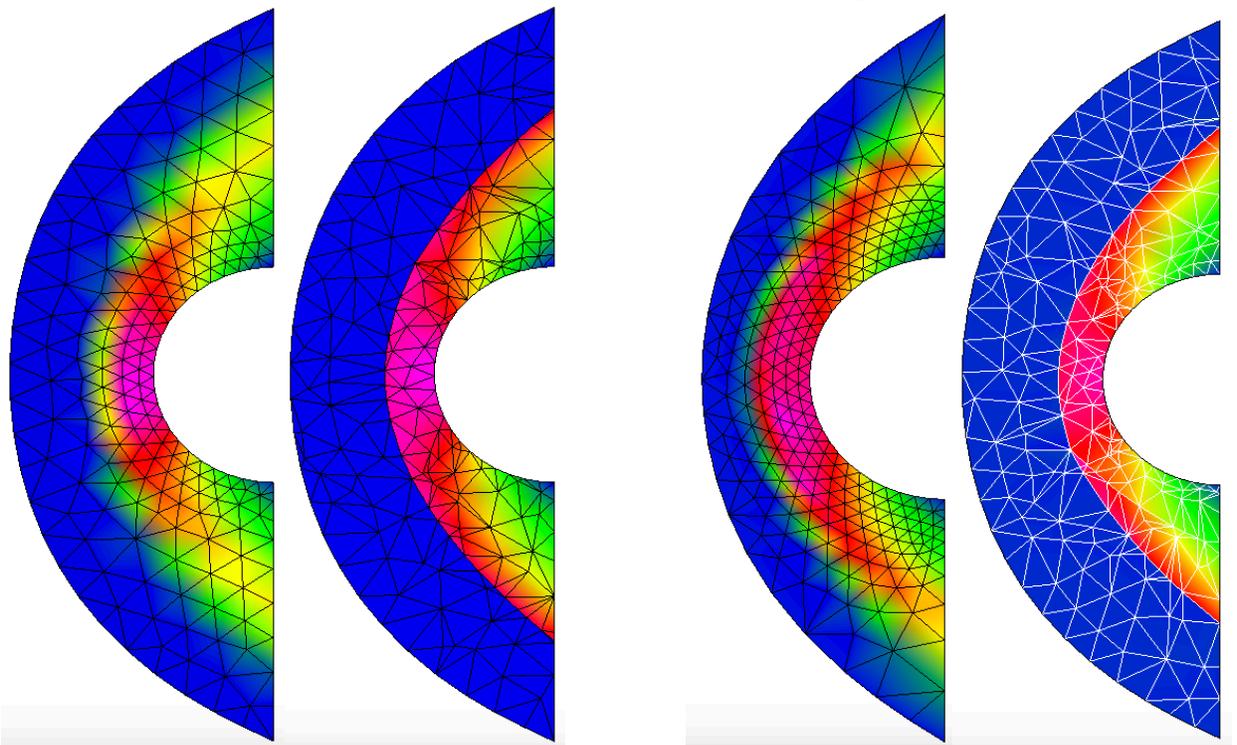

**Figure 28. Initial grid and density field obtained by the DG(P0) solution and final converged mesh and density field obtained by the MDG(P1) solution for a Mach number 8 (left) and 18 (right) hypersonic flow**

## V. Concluding Remarks and Future Work

A MDG+ICE method based on a continuous variational formulation for the interface conservation and a space-time discontinuous Galerkin formulation for the conservative laws has been developed for solving the compressible Euler equations. A number of test cases for both 1D unsteady and 2D steady compressible Euler equations have been conducted to assess the accuracy and robustness of the MDG+ICE method. The preliminary results are highly promising and encouraging, demonstrating that an exponential rate of convergence for Sod and Lax-Harden shock tube problems can be achieved and highly accurate solutions without overheating to both double-rarefaction wave and Noh problems can be obtained, something that no other numerical methods can do to the best of our knowledge. Indeed, by recognizing the importance of the interface conservation and by enforcing the interface conservation via grid movement and grid management, the MDG+ICE method, 1. can automatically detect and fit all types of discontinuities virtually exactly and resolve accurately solutions with discontinuous derivatives, which in turn allows the MDG+ICE method to achieve the designed optimal rate of convergence even for discontinuous solutions and solutions with singularities; 2. does not require any exact or approximate Riemann solver-based numerical flux function to achieve stability nor does it require any strategies (limiter/ENO-WENO/artificial viscosity) to suppress spurious oscillations in the vicinity of strong discontinuities; and 3. is extremely robust in the sense that a). It can handle highly distorted and even tangled meshes with zero or even negative volumes; and b). Even negative pressure and negative density do not lead to the breakdown of a solution process. Future work will be focused on exploring and developing numerical algorithms to solve the resultant nonlinear least-squares problems efficiently, effectively, and reliably. The extension of the MDG+ICE method for solving viscous flow, reactive flow, and multi-material flow problems is also underway.

## Acknowledgement

This work was performed under the auspices of the U.S. Department of Energy by Lawrence Livermore National Laboratory under Contract DE-AC52-07NA27344, and funded by the Laboratory Directed Research and Development Program at LLNL under project tracking code 19-ERD-015. We would like to thank Dr. Andrew Corrigan at Naval Research Laboratory for the fruitful discussion.## References

[1] W.H. Reed and T.R. Hill, Triangular Mesh Methods for the Neutron Transport Equation, **Los Alamos Scientific Laboratory Report**, LA-UR-73-479, 1973.

[2] B. Cockburn, S. Hou, and C. W. Shu, TVD Runge-Kutta Local Projection Discontinuous Galerkin Finite Element Method for conservation laws IV: the Multidimensional Case, **Mathematics of Computation**, Vol. 55, pp. 545-581, 1990.eorge, P. L., Automatic Mesh Generation, J. Wiley & Sons, 1991.

[3] B. Cockburn, and C. W. Shu, The Runge-Kutta Discontinuous Galerkin Method for conservation laws V: Multidimensional System, **Journal of Computational Physics**, Vol. 141, pp. 199-224, 1998.

[4] B. Cockburn, G. Karniadakis, and C. W. Shu, The Development of Discontinuous Galerkin Method, in Discontinuous Galerkin Methods, Theory, Computation, and Applications, edited by B. Cockburn, G.E. Karniadakis, and C. W. Shu, Lecture Notes in Computational Science and Engineering, Springer-Verlag, New York, 2000, Vol. 11 pp. 5-50, 2000.

[5] F. Bassi and S. Rebay, High-Order Accurate Discontinuous Finite Element Solution of the 2D Euler Equations, **Journal of Computational Physics**, Vol. 138, pp. 251-285, 1997.

[6] H. L. Atkins and C. W. Shu, Quadrature Free Implementation of Discontinuous Galerkin Method for Hyperbolic Equations, **AIAA Journal**, Vol. 36, No. 5, 1998.

[7] T. C. Warburton, and G. E. Karniadakis, A Discontinuous Galerkin Method for the Viscous MHD Equations, **Journal of Computational Physics**, Vol. 152, pp. 608-641, 1999.

[8] J. S. Hesthaven and T. Warburton, Nodal Discontinuous Galerkin Methods: Algorithms, Analysis, and Applications, Texts in **Applied Mathematics**, Vol. 56, 2008.

[9] A. Pandare and H. Luo, A Hybrid Reconstructed Discontinuous Galerkin and Continuous Galerkin method for Incompressible Flows on Unstructured Grids, **Journal of Computational Physics**, Vol. 322, pp. 491-510, 10.1016/j.jcp.2016.07.002, 2016.

[10] B. Karami Halashi, and H. Luo, A Discontinuous Galerkin Method for Magnetohydrodynamics on Arbitrary Grids, **Journal of Computational Physics**, Vol. 326, pp. 258-277, 10.1016/j.jcp.2016.08.055, 2016.

[11] K. J. Fidkowski, T. A. Oliver, J. Lu, and D. L. Darmofal, $p$-Multigrid solution of high-order discontinuous Galerkin discretizations of the compressible Navier–Stokes equations, **Journal of Computational Physics**, Vol. 207, No. 1, pp. 92-113, 2005.

[12] H. Luo, J. D. Baum, and R. Löhner, A Discontinuous Galerkin Method Using Taylor Basis for Compressible Flows on Arbitrary Grids, **Journal of Computational Physics**, Vol. 227, No 20, pp. 8875-8893, October 2008.

[13] H. Luo, J.D. Baum, and R. Löhner, On the Computation of Steady-State Compressible Flows Using a Discontinuous Galerkin Method, **International Journal for Numerical Methods in Engineering,** Vol. 73, No. 5, pp. 597-623, 2008.27